%% file: stability.tex
\documentclass[11pt,a4paper,bibtotoc]{scrartcl}
\usepackage[english]{babel}
\usepackage{amsmath}
\usepackage{amsthm}
\usepackage{amssymb}
\usepackage{graphics}
\usepackage{alltt}
\usepackage{latexsym}
\usepackage{pstricks}
\usepackage{graphicx}
\usepackage{algorithm}
\usepackage{algorithmic}
\usepackage{subfigure}
\usepackage{multirow}
\usepackage{multicol}

\def\N{\mathbb{N}}
\def\R{\mathbb{R}}
\def\Z{\mathbb{Z}}
\def\C{\mathbb{C}}
\def\T{\mathbb{T}}

\def\ti{\mbox{\scriptsize{\rm i}}}

\newcommand{\eip}[1]{{\rm e}^{ 2\pi{\ti} #1}}
\newcommand{\eim}[1]{{\rm e}^{-2\pi{\ti} #1}}
\newcommand{\zb}[1]{\ensuremath{\boldsymbol{#1}}}
\newcommand{\adj}{{\vdash \hspace*{-1.72mm} \dashv}}

\newcommand{\myspan}{{\rm span}}
\newcommand{\dist}{{\rm dist}}
\newcommand{\trace}{{\rm trace}}
\newcommand{\cond}{{\rm cond}}
\newcommand{\diag}{{\rm diag}}

\renewcommand{\d}{{\rm d}}
\newcommand{\supp}{{\rm supp}}

\newcommand{\indexset}{{I_N}}

\renewcommand{\Box}{\hspace*{0ex} \hfill \rule{1.5ex}{1.5ex} \\ \goodbreak}
\newcommand{\Boxgl}{\par\vspace{-5ex} \hspace*{0ex} \hfill
  \rule{1.5ex}{1.5ex} \\ \goodbreak\goodbreak}
\newcommand{\bend}{\hspace*{0ex} \hfill \hbox{\vrule height
    1.5ex\vbox{\hrule width 1.4ex \vskip 1.4ex\hrule  width 1.4ex}\vrule
    height 1.5ex} \goodbreak}
\newtheorem{theorem}{Theorem}[section]
\newtheorem{lemma}[theorem]{Lemma}
\newtheorem{remark}[theorem]{Remark}
\newtheorem{definition}[theorem]{Definition}
\newtheorem{example}[theorem]{Example}
\newtheorem{corollary}[theorem]{Corollary}

\newenvironment{Theorem}{\goodbreak \begin{theorem}\sl}{\end{theorem}}
\newenvironment{Lemma}{\goodbreak \begin{lemma}\sl}{\end{lemma}}
\newenvironment{Remark}{\goodbreak \begin{remark}\rm}{\bend\end{remark}}

\newenvironment{Definition}{\goodbreak \begin{definition}\rm}{\bend\end{definition}}
\newenvironment{Corollary}{\goodbreak \begin{corollary}\rm}{\end{corollary}}

\numberwithin{equation}{section}
\numberwithin{table}{section}
\numberwithin{figure}{section}

\title{Stability Results for Scattered Data Interpolation by Trigonometric
  Polynomials}

\author{Stefan Kunis and Daniel Potts
  \thanks{Department of Mathematics, Chemnitz University of Technology,
    09107 Chemnitz, Germany (\{kunis,potts\}@mathematik.tu-chemnitz.de)}}

\date{}

\begin{document}
\maketitle

\begin{abstract}
A fast and reliable algorithm for the optimal interpolation of scattered data
on the torus $\T^d$ by multivariate trigonometric polynomials is presented.
The algorithm is based on a variant of the conjugate gradient method in
combination with the fast Fourier transforms for nonequispaced nodes.
The main result is that under mild assumptions the total complexity for
solving the interpolation problem at $M$ arbitrary nodes is of order ${\cal
  O}(M\log M)$.
This result is obtained by the use of localised trigonometric kernels
where the localisation is chosen in accordance to the spatial dimension $d$.
Numerical examples show the efficiency of the new algorithm.

\smallskip


\end{abstract}

\section{Introduction}
We discuss the approximation of scattered data by $d$-variate periodic
functions $f:\T^d\rightarrow\C$, where $\T:=[-\frac{1}{2},\frac{1}{2})$
denotes the torus.
In practical applications we are often confronted with the situation that
experimental data or measured values of a function are only known at a finite sampling
set ${\cal X} := \{\zb x_j \in \T^d \,:\, j=0,\hdots,M-1\}$.
Especially, nonuniform sampling sets appear in more and more applications in
recent years.
Given a notion of the distance of two points by $\dist \left(\zb x,\zb
  x_0\right) := \min_{\zb j \in \mathbb Z^d} \left\| (\zb x+\zb j) - \zb
  x_0\right\|_{\infty}$, we measure the ``nonuniformity'' of $\cal X$ by the
mesh norm and the separation distance, defined by
\begin{equation*}
  \delta:=2\max_{\zb x\in\T^d}\min_{j=0,\hdots,M-1}\dist(\zb x_j,\zb x),\qquad
  q:=\min\limits_{j,l=0,\hdots,M-1;j\ne l}\dist\left(\zb x_j,\zb x_l\right),
\end{equation*}
respectively.
Obviously, the relation $q\le M^{-1/d}\le\delta$ is fulfilled.

For given samples $(\zb x_j,y_j)\in \T^d\times\C,\; j=0,\hdots,M-1$, a
polynomial degree $N\in 2\N$, and the index set
$\indexset:=\left\{-\frac{N}{2},\hdots,\frac{N}{2}-1\right\}^d$ of
frequencies, we construct a $d$-variate trigonometric polynomial
\begin{equation*}
  f\left(\zb x\right):=\sum\limits_{\zb k \in \indexset} \hat f_{\zb k}
  \eip{\zb k \zb x} 
\end{equation*}
such that  $f(\zb x_j) \approx y_j,\;j=0,\hdots,M-1$.
Turning this into matrix vector notation, we aim to solve the system of linear
equations
\begin{equation}\label{eq:1}
 \zb A \zb {\hat f} \approx \zb y
\end{equation}
for the unknown vector of Fourier coefficients $\zb {\hat f}:=(\hat f_{\zb k}
)_{\zb k\in \indexset} \in \C^{N^d}$.
Throughout the paper, we denote the vector of the given sample values by $\zb
y:=(y_{j})_{j=0,\ldots,M-1}\in \mathbb{C}^{M}$ and the nonequispaced Fourier
matrix by
\begin{equation*}
  \zb A=\zb A_{\cal X}:=\left(\eip{\zb k \zb x_j}\right)_{j=0,\hdots,M-1;\zb k
  \in \indexset} \in \C^{M\times N^d}.
\end{equation*}

In contrast to the widely used nonequispaced FFT for the fast matrix vector
multiplication with $\zb A$, see \cite{postta01} and its references, the
{\em efficient} solution of \eqref{eq:1} is still a challenging goal.
Besides recently developed sparse reconstruction techniques, see e.g.
\cite{DaDeDeM04,KuRa06} and their references, a standard method to determine
$\zb {\hat f}$ is to solve the general linear least squares problem 
$\|\zb {\hat f}\|_2 \rightarrow \min$ subject to $\|\zb y - \zb A \zb
{\hat f}\|_2=\min$, see, e.g., \cite[p. 15]{Bj96}.
This can be done by means of the singular value decomposition which is not
practical in the present situation for large problems due to its time and
memory requirements.
Direct solvers for the univariate case $d=1$ in \cite{fas97,ReAmGr} obtain a
solution in ${\cal O}(NM)$ floating point operations.

For $N^d<M$, the linear system \eqref{eq:1} is over-determined, so that in
general the given data $\zb y$ will be only approximated up to a
{\em residual} $\zb r:=\zb y-\zb A \zb {\hat f}$. 
In order to compensate for clusters in the sampling set ${\cal X}$, it is
useful to incorporate weights $w_j> 0$ and to consider the weighted
approximation problem
\begin{equation} \label{eq:ap}
  \|\zb y - \zb A \zb {\hat f}\|_{\zb W}^2
  = \sum_{j=0}^{M-1} w_j |y_j-f(\zb x_j)|^2
  \stackrel{\zb {\hat f}}{\rightarrow} \min,
\end{equation}
where $\zb W:=\diag(w_j)_{j=0,\ldots,M-1}$.
In \cite{Groechenig92} it has been proven that this problem has a unique
solution if $N<(\frac{\pi}{\log 2} \,d \,\delta)^{-1}$.
Its solution is computed iteratively by means of the conjugate gradient method
in \cite{FeGrSt95,BaGr03,GrSt04}, where the multilevel Toeplitz structure of
$\zb A^{\adj}\zb W\zb A$ is used for fast matrix vector multiplications.
Slightly more stable with respect to rounding errors is the CGNR method, cf.
\cite[pp. 288]{Bj96}, which iterates the original residual $\zb r_l=\zb y -
\zb A \zb {\hat f}_l$ instead of the residual $\zb A^{\adj}\zb W \zb r_l$ of
the normal equations.
Note furthermore, that it has been suggested in \cite{RaSt98} to incorporate
some ``knowledge on the decay of the Fourier coefficients''
$\zb{\hat W}:= \diag(\hat w_{\zb k})_{\zb k\in \indexset}$,
$\hat w_{\zb k}> 0$.
Their approach is based on the weighted least squares problem
\begin{equation*}
  \left\|\zb A^{\adj} \zb W \left(\zb y - \zb A \zb {\hat
    f}\right)\right\|_{\zb {\hat W}^{-1}} \stackrel{\zb {\hat
      f}}{\rightarrow} \min.
\end{equation*}

In contrast, we focus on the under-determined and consistent linear system
$\zb A \zb {\hat f} = \zb y$, i.e., we expect to interpolate the given data
$y_j\in\C$, $j=0, \hdots,M-1$, exactly.
We show that the nonequispaced Fourier matrix $\zb A$ has full rank $M$ for
every polynomial degree $N>2 \,d \,q^{-1}$.
In particular, we incorporate {\em damping factors} $\hat w_{\zb k}> 0$,
$\zb k \in I_N$, and consider the {\em optimal interpolation problem}
\begin{equation} \label{eq:optInterp}
  \|\zb {\hat f}\|_{\zb {\hat W^{-1}}}^2
  =\sum_{\zb k\in \indexset} \frac{|\hat f_{\zb k}|^2}{\hat w_{\zb k}}
  \stackrel{\zb {\hat
      f}}{\rightarrow} \min \quad \text{subject to} \quad \zb A \zb {\hat
    f} = \zb y,
\end{equation}
where $\zb{\hat W}:= \diag(\hat w_{\zb k})_{\zb k\in \indexset}$.
We prove that for a large class of ``smooth'' damping factors $\hat w_{\zb k}$
problem \eqref{eq:optInterp} is well-conditioned, where the ``smoothness'' has
to be chosen with respect to the spatial dimension $d$.
We propose to solve problem \eqref{eq:optInterp} by a version of the
conjugate gradient method in combination with the nonequispaced FFT
\cite{st97,postta01,kupo02C} to efficiently perform each iteration step.

The outline of this paper is as follows:
In Section \ref{sect:optInterp} we set up the basic notation and relate the
optimal interpolation problem \eqref{eq:optInterp} to a particular
trigonometric kernel.
Furthermore, we propose Algorithm \ref{algo:CGNE} for computing the solution
of the interpolation problem efficiently.
For the sake of analysing the convergence of this algorithm, we then present
our theory on localised trigonometric kernels in Section 3.
Our first result in Theorem~\ref{theorem:T1_smooth} is a version of the
typical smoothness-decay principle in Fourier analysis and relates the
``smoothness'' of the weights in \eqref{eq:optInterp} to the localisation of
the corresponding trigonometric kernel. 
We use this decay in Section 4 to prove that well separated sampling nodes
yield a stable interpolation problem \eqref{eq:optInterp}.
The eigenvalue estimates are given for the univariate setting in Theorem
\ref{theorem:T1allgemein} and for the multivariate setting in Theorem
\ref{theorem:Tdallgemein}.
Subsequently, Corollary \ref{cor:qsep} applies the general result to a
particular class of damping factors and concludes with conditions sufficient
for the full rank of $\zb A$.
As the equidistant case in Theorem \ref{theorem:eq} and Corollary
\ref{cor:eq} reveals, the assumption on the separation distance is of optimal
order.
We provide numerical examples in Section 5 and draw our conclusion in
Section 6.
The software and all numerical examples are available from our NFFT-homepage
\cite{kupo02C}.

\section{Optimal interpolation and its iterative solution}
\label{sect:optInterp}
After setting up our notation in Definition \ref{def:kernel}, we prove in
Lemma \ref{lemma:ne} that the optimal interpolation problem
\eqref{eq:optInterp} can be stated as normal equations and the matrix in
these equations obeys special structure.
Furthermore, we propose Algorithm \ref{algo:CGNE} for the iterative solution
of the interpolation problem and state a basic convergence result for this
scheme.

\begin{Definition}  \label{def:kernel}
  Let $d\in\N,\;N\in 2\N$, and $I_N=\{-\frac{N}{2},\hdots,\frac{N}{2}-1\}^d$
  be given. 
  We define for positive weights $\hat w_{\zb k}>0,\; \zb k \in
  \indexset$, with normalisation $\sum_{\zb k\in I_N} \hat w_{\zb k}=1$ and for $\zb x
  \in \T^d$ the trigonometric {\em kernel}
  \begin{equation*}
    K_N \left(\zb x\right):=\sum\limits_{\zb k \in \indexset} \hat w_{\zb k}
    \eip{\zb k \zb x}.
  \end{equation*}
  The particular class of {\em tensor product kernels} is given by
  \begin{equation*}
    K_N\left(\zb x\right)=\prod\limits_{t=0}^{d-1}
    \tilde K_N\left(x_t\right)
  \end{equation*}
  where $\tilde K_N$ denotes a univariate kernel and $\zb
  x=(x_0,\hdots,x_{d-1})^{\top}$.

  Furthermore, given a sampling set ${\cal X}\subset\T^d$, we define the
  {\em kernel matrix}
  \begin{equation}\label{def:K}
    \zb K_N: = \left(K_N (\zb x_j-\zb x_l) \right)_{j,l=0,\ldots,M-1}\in
    \C^{M\times M}.
  \end{equation}
  We denote by $\Lambda=\Lambda\left(\zb K_N\right)$ and
  $\lambda=\lambda\left(\zb K_N\right)$ the largest and smallest eigenvalue of 
  the kernel matrix $\zb K_N$, respectively.
  Their ratio is denoted by the condition number $\cond(\zb K_N)=
  \frac{\Lambda}{\lambda}$.
\end{Definition}

Note, that from the definition immediately follows $K_N(\zb 0)=\max_{\zb
  x\in\T^d} |K_N(\zb x)|=1$ and $(\zb K_N)_{j,j}=1$, $j=0,\hdots,M-1$.
The following theorem collects some basic facts.
\begin{Lemma}\label{lemma:ne}
  Let the number of nodes $M\in\N$, the sampling set ${\cal X}\subset\T^d$,
  the polynomial degree $N \in 2\N$, and the damping factors $\hat w_{\zb
  k}>0,\; \zb k \in \indexset$, be given.
  The optimal interpolation problem \eqref{eq:optInterp} is equivalent
  to the {\em damped normal equations of second kind}
  \begin{equation}\label{eq:ne}
    \zb K_N \zb {\tilde f} = \zb y,\qquad \zb {\hat
      f}=\zb {\hat W} \zb A^{\adj} \zb {\tilde f},
  \end{equation}
  where the kernel matrix $\zb K_N\in \C^{M\times M}$ obeys the factorisation
  \begin{equation}\label{eq:fact}
    \zb K_N = \zb A \zb {\hat W} \zb A^{\adj},
  \end{equation}
  hence is positive semidefinite.
\end{Lemma}
Proof.
The second assertion follows from $(\zb A \zb
{\hat W} \zb A^{\adj})_{j,l}=\sum_{\zb k \in \indexset} \eip{\zb k \zb x_j}
\hat w_{\zb k} \eim{\zb k \zb x_l}$ and Definition \ref{def:kernel}.
Furthermore, a solution $\zb {\hat f}$ of $\zb A \zb {\hat f}=\zb y$ has
minimal weighted norm if and only if it is perpendicular with respect to the
weights to the null-space of $\zb A$, i.e., $\zb {\hat W}^{-{1}/{2}} \zb {\hat
  f} \perp {\cal N}(\zb A \zb {\hat W}^{{1}/{2}})$.
We conclude \eqref{eq:ne} by the fact that the orthogonal complement of
the null-space of a matrix is just the range of its adjoint.
 \Box

Denoted in Algorithm \ref{algo:CGNE} by CGNE, cf. \cite[pp. 288]{Bj96}, we
solve the {\bf N}ormal equations \eqref{eq:ne} by the {\bf C}onjugate {\bf
  G}radient method, minimising in each iteration the {\bf E}rror.

\begin{algorithm}[h!]
  \caption{CGNE}\label{algo:CGNE}
  \begin{tabular}{ll}
    Input: & dimension $d\in \N$, number of samples $M\in \N$, polynomial
    degree $N\in 2\N$;\\
    & sampling set ${\cal X}\subset\T^d$, samples $\zb y \in \C^M$, and initial vector $\zb
    {\hat f}_0 \in \C^{N^d}$
  \end{tabular}
  \begin{algorithmic}
    \topsep=1.5ex \itemsep=1.5ex
    \STATE
    \STATE $\zb {r}_0=\zb y - \zb A \zb {\hat f}_0$
    \STATE $\zb {\hat p}_0=\zb A^{\adj} \zb {r}_0$
    \FOR{$l=0,\hdots$}
    \topsep=1.5ex \itemsep=1.5ex
    \STATE $\alpha_l= \zb {r}_l^{\adj} \zb {r}_l \,/\,
    {\zb {\hat p}_l^{\adj} \zb {\hat W} \zb {\hat p}_l}$
    \STATE $\zb {\hat f}_{l+1}=\zb {\hat f}_{l}
    +\alpha_l \zb {\hat W} \zb {\hat p}_l$
    \STATE $\zb {r}_{l+1}=\zb {r}_{l}- \alpha_l \zb A \zb {\hat W} \zb {\hat
    p}_l$ 
    \STATE $\beta_l={\zb {r}_{l+1}^{\adj} \zb {r}_{l+1}}\,/\,
    {\zb {r}_l^{\adj} \zb {r}_l}$
    \STATE $\zb {\hat p}_{l+1}=\beta_l \zb {\hat p}_l + \zb A^{\adj} \zb
    {r}_{l+1}$ 
    \ENDFOR
    \STATE
  \end{algorithmic}
  \begin{tabular}{ll}
    Output: & the $l$-th iterate $\zb {\hat f}_{l}$
  \end{tabular}
\end{algorithm}

The proposed method finds approximations from a Krylov space closely related
to the one of the CGNR method for \eqref{eq:ap}, but with minimal error
instead of minimal residual.
Note that we exploit the factorisation in \eqref{eq:fact} to iterate the
original vector $\zb {\hat f}$ instead of the vector $\zb {\tilde f}$,
cf. equation \eqref{eq:ne}.
Hence, we use fast matrix vector multiplications for $\zb A$ and $\zb
A^{\adj}$ by means of the fast Fourier transforms at nonequispaced nodes
(NFFT) having an arithmetical complexity of ${\cal O}(N^d\log(N^d)+M |\log \epsilon|^d )$ in
 each iteration, where $\epsilon$ is the prescribed accuracy.
Details concerning NFFT algorithms can be found for example in
\cite{st97,postta01} and a corresponding software package in \cite{kupo02C}.
Applying the standard estimate for the convergence of the conjugate gradient
method we obtain the following lemma.

\begin{Lemma}\label{lemma:cg_error}
Let the kernel matrix $\zb K_N$ in \eqref{def:K} be regular and let
 $\zb {\hat e}_l:=\zb {\hat f}_l-\zb {\hat W}\zb A^{\adj}\zb K_N^{-1} \zb y$
 denote the error of the $l$-th iterate within Algorithm \ref{algo:CGNE}.
Then the a-priori error bound
\begin{equation*}
  \left\|\zb {\hat e}_l\right\|_{\zb {\hat W}^{-1}} \le
  2 \left(\frac{\sqrt{\Lambda}-\sqrt{\lambda}}{\sqrt{\Lambda}+
  \sqrt{\lambda}}\right)^l 
  \left\|\zb {\hat e}_0\right\|_{\zb {\hat W}^{-1}}
\end{equation*}
holds true.
\end{Lemma}
Proof. We note that $\|\zb {\hat e}_l\|_{\zb {\hat W}^{-1}}=\|\zb {\tilde f}_l
- \zb K_N^{-1} \zb y\|_{\zb K_N}$, where $\zb {\tilde f}_l$ denotes the $l$-th
iterate of the conjugate gradient method applied to equation $\zb K_N \zb
{\tilde f} = \zb y$, cf. Lemma \ref{lemma:ne}, and apply the standard
estimate for the conjugate gradient method, see also \cite[pp. 288]{Bj96}. \Box

This result includes the special case of $M=N^d$ equispaced nodes and no
damping, i.e. $\hat w_{\zb k}=1$, $\zb k \in I_N$, where the first iterate of
our algorithm is already the solution to equation \eqref{eq:ne}.
We present estimates for the extremal eigenvalues $\lambda,\Lambda$ dependent
only on the quantities $N,q$, and the damping factors $\hat w_{\zb k},\; \zb k
\in \indexset$.
Analogous results for the stability of the interpolation by radial and zonal
functions are obtained in \cite{NaSiWa98,Wen05}.
Section \ref{sect:kernels} prepares our estimates by constructing localised
kernels.

\begin{Remark}\label{remark:ref}
Before that, we would like to comment on the following:

The weighted norm in \eqref{eq:optInterp} is induced by the inner product
$\zb {\hat g}^{\adj} \zb {\hat W}^{-1} \zb {\hat f}$.
In particular, the definition $
  \left\langle f,g\right\rangle_{\zb {\hat W}^{-1}}:=
  \zb {\hat g}^{\adj} \zb {\hat W}^{-1} \zb {\hat f}$
makes the space of trigonometric polynomials
$T_N:=\myspan\left\{\eip{\zb k \cdot}:\;\zb k \in \indexset\right\}$ to a
reproducing kernel Hilbert space.
Its reproducing kernel is given by $K_N$, i.e., the point evaluations obey
$f(\zb x) = \langle f,K_N(\cdot - \zb x)\rangle_{\zb {\hat W}^{-1}}$.

Moreover, the solution $f(\zb x)=\sum_{\zb k \in \indexset} \hat f_{\zb k}
\eip{\zb k \zb x}$ of the normal equations \eqref{eq:ne} has comparable norm
to the given samples, i.e.,
\begin{equation*}
  \Lambda^{-1}\left\|\zb y\right\|_2^2 \le
  \left\langle f,g\right\rangle_{\zb {\hat W}^{-1}} \le
  \lambda^{-1}\left\|\zb y\right\|_2^2.
\end{equation*}
This norm equivalence is due to fact, that the field of values of the matrix
$\zb K_N^{-1}$ is bounded by its extremal eigenvalues and 
$
  \zb y^{\adj} \zb K_N^{-1} \zb y
  = \zb {\tilde f}^{\adj} \zb K_N \zb {\tilde f}
  = \zb{\hat f}^{\adj} \zb {\hat W}^{-1} \zb{\hat f}
$.
\end{Remark}

\section{Localised kernels}\label{sect:kernels}
Starting from a class of admissible weight functions in Definition
\ref{def:BV}, we construct localised trigonometric kernels in Theorem
\ref{theorem:T1_smooth}, where Lemma \ref{lemma:decay_g} serves as an
intermediate step.
Following the smoothness-decay principle in Fourier analysis, we relate the
smoothness of the weight function to the decay of the kernel $K_N$
built upon the sampled weights.
A related approach is taken in \cite[Thm. 2.2]{MhPr00} for the detection
of singularities.
The particular class of B-Spline kernels, cf. Definition \ref{def:spline}, is
considered in Corollary \ref{cor:spline}.
While we present our results on the connection between smooth weight functions
and localised kernels for the univariate case, we give its generalisation to
the class of tensor product kernels in Corollary \ref{cor:tensorproduct}.

\begin{Definition}\label{def:BV}
  For $\beta \in \N$, $\beta\ge 2$, a continuous function $g:\R\rightarrow \R$
  is an {\em admissible weight function of order $\beta$} if it is
  nonnegative, possesses a $(\beta-1)$-fold derivative $g^{(\beta-1)}$
  of bounded variation, i.e.,
  \begin{equation*}
    \left|g^{\left(\beta-1\right)}\,\right|_V
    := \int\limits_{\R}
    \left|\d g^{\left(\beta-1\right)}\left(z\right) \right |=
    \sup
    \sum\limits_{j=0}^{n-1} \left|
    g^{\left(\beta-1\right)}\left(z_{j+1}\right) -
    g^{\left(\beta-1\right)}\left(z_{j}\right) \right| 
    < \infty,
  \end{equation*}
  where the supremum is taken over all strictly increasing real sequences
  $\{z_j\}_{j\in\N_0}$, and satisfies the additional properties $\supp\, g=
  [-\frac{1}{2},\frac{1}{2}]$, $g^{(\gamma)}(\pm \frac{1}{2})=0$ for
  $\gamma=0,\ldots,\beta-1$, $g(z)>0,\,|z|<\frac{1}{2}$, and the
  normalisation $\|g\|_{L^1}=1$.
  We denote by $BV^{\beta-1}_0$ the set of admissible weight functions of
  order $\beta$.

  Furthermore, we define for notational convenience the {\em zeta function}
  $\zeta(\beta):=\sum_{r=1}^{\infty} r^{-\beta}$, $\beta>1$, and for $g\in
  BV_0^{\beta-1}$ the {\em norm of the samples} 
  \begin{equation*}
    \left\|g\right\|_{1,N}:=\sum_{k=-\frac{N}{2}}^{\frac{N}{2}}
    g\left(\frac{k}{N}\right)\,.
  \end{equation*}
\end{Definition}

The following lemma prepares Theorem \ref{theorem:T1_smooth}.

\begin{Lemma}\label{lemma:decay_g}
  For $\beta \in \mathbb{N},\;\beta\ge 2$, let a function $g\in BV^{\beta-1}_0$
  be given.
  Then for $N\in 2\N$, $N\ge 2\beta$, and $x \in \left[-\frac{1}{2},
  \frac{1}{2}\right]\setminus\left\{0\right\}$ the following estimates
  hold true
  \begin{eqnarray*}
    \left| \sum\limits_{k=-\frac{N}{2}}^{\frac{N}{2}}
    g\left(\frac{k}{N}\right) \eip{kx} \right|
    & \le &
    \frac{\left(2^{\beta}-1\right)\zeta\left(\beta\right)
    \left|g^{\left(\beta-1\right)}\right|_V}  
    {\left(2N\right)^{\beta-1} |2\pi x|^{\beta}},\\
    \left\|g\right\|_{1,N}
    & \ge & 
    N\left(1-2\zeta\left(\beta\right) \left(4\pi \beta\right)^{-\beta}
    \left|g^{\left(\beta-1\right)}\right|_V\right). 
  \end{eqnarray*}
\end{Lemma}

Proof. First, we define for $x,z \in \left[-\frac{1}{2},\frac{1}{2}\right]$ the
function $h_x\left(z\right):=  g\left(z\right)\eip{N x z}$.
Thus, the Poisson summation formula yields
\begin{equation*}
  \frac{1}{N}
  \sum\limits_{k=-\frac{N}{2}}^{\frac{N}{2}}
  g\left(\frac{k}{N}\right) \eip{kx} =
  \frac{1}{N}
  \sum\limits_{k=-\frac{N}{2}}^{\frac{N}{2}}
  h_x\left(\frac{k}{N}\right)\, =\, \sum_{r\in \mathbb{Z}}
  \int\limits_{-\frac{1}{2}}^{\frac{1}{2}}
  h_x\left(z\right)  \eim{N r z} \d z
\end{equation*}
and by applying integration by parts and the fact that
$g^{\left(\gamma\right)}\left(\pm\frac{1}{2}\right)=0$ for $\gamma=
0,\ldots,\beta-2$ further
\begin{eqnarray*}
  \left| \sum\limits_{k=-\frac{N}{2}}^{\frac{N}{2}}
  h_x\left(\frac{k}{N}\right) \right|
  &=&
  \left| N \sum_{r\in \mathbb{Z}} \left(2\pi {\rm i} N
      \left(r-x\right)\right)^{-\left(\beta-1\right)} 
    \int\limits_{-\frac{1}{2}}^{\frac{1}{2}}
    g^{\left(\beta-1\right)}\left(z\right)  \eip{N z \left(x - r\right)}
    \d z \right |\\
  &=&
  \left| \frac{N}{\left(2\pi {\rm i} N\right)^{\beta}}
  \sum_{r\in \mathbb{Z}} \left(r-x\right)^{-\beta} 
  \int\limits_{-\frac{1}{2}}^{\frac{1}{2}}
    g^{\left(\beta-1\right)}\left(z\right) \left(\frac{\d}{\d z}\,
    \eip{N z \left(x - r\right)}\right)  \d z \right |\\
  &\le&
  \frac{1+\left|x\right|^{\beta} \sum\limits_{r\in \Z\setminus\{0\}}
    \left|r-x\right|^{-\beta}}
  {\left(2\pi\right)^{\beta} N^{\beta-1} |x|^{\beta}} \sup_{r_0\in\Z}
  \left | \, \int\limits_{-\frac{1}{2}}^{\frac{1}{2}}
    g^{\left(\beta-1\right)}\left(z\right)   \left(\frac{\d}{\d z}\,\eip{N z \left(x -
          r_0\right)} \right) 
    \d z \right |\,.
\end{eqnarray*}
Using $1+\left|x\right|^{\beta} \sum_{r\in \Z\setminus\{0\}}
\left|r-x\right|^{-\beta} \le (2^{\beta}-1)2^{1-\beta}\zeta(\beta)$ for
$|x|\le\frac{1}{2}$ and
\begin{equation*}
  \left |
   \int\limits_{-\frac{1}{2}}^{\frac{1}{2}}
   g^{\left(\beta-1\right)}\left(z\right)   \left(\frac{\d}{\d z}\,
   \eip{N z \left(x - r_0\right)}\right) \d z \right |
   \le
   \left|g^{\left(\beta-1\right)}\,\right|_V
\end{equation*}
yields the assertion.

By the Poisson summation formula, we note furthermore that
\begin{equation*}
  \frac{1}{N} \left\|g\right\|_{1,N}
  \ge 1-\left|\sum\limits_{r\in\mathbb{Z}\setminus\left\{0\right\}}
  \int\limits_{-\frac{1}{2}}^{\frac{1}{2}} g\left(z\right) \eim{Nrz}\d
  z\right|
\end{equation*}
and proceed analogously in order to prove the second assertion where we use
$N\ge 2\beta$ to obtain an estimate independent of $N$. \Box

\begin{Theorem}
 \label{theorem:T1_smooth}
 For $\beta \in \mathbb{N},\;\beta\ge 2$, let a function $g\in BV^{\beta-1}_0$
 be given.
 Furthermore, let $N\in2\mathbb{N}$, $N\ge 2\beta$, and the 
 damping factors
 \begin{equation*}
   \hat w_k=\frac{g\left(\frac{k}{N}\right)+g\left(\frac{k+1}{N}\right)}{
     2\left\|g\right\|_{1,N}}\,,\qquad k=-\frac{N}{2},\hdots,\frac{N}{2}-1,
 \end{equation*}
 be given.
 Then the kernel $K_N$, cf. Definition \ref{def:kernel},
 fulfils
 \begin{equation*}
   \left| K_N\left(x\right)\right| \le 
   \frac{\left(2^{\beta}-1\right)\zeta\left(\beta\right)
     \left|g^{\left(\beta-1\right)}\right|_V}
   {2^{\beta-1}\left(2\pi\right)^{\beta}- \zeta\left(\beta\right)
     \beta^{-\beta} \left|g^{\left(\beta-1\right)}\right|_V}
   \frac{1}{ N^{\beta}
     |x|^{\beta}}
 \end{equation*}
 for $x\in\left[-\frac{1}{2},\frac{1}{2}\right]\setminus
 \left\{0\right\}$.
\end{Theorem}
Proof. Note first, that
\begin{equation*}
K_N\left(x\right)
= \frac{1+\eim{x}}{2\left\|g\right\|_{1,N}}
\sum\limits_{k=-\frac{N}{2}}^{\frac{N}{2}} g\left(\frac{k}{N}\right) \eip{kx}.
\end{equation*}
Thus, we obtain $K_N(0)=1$ and by applying Lemma \ref{lemma:decay_g} also
the decay property.
\Box

We apply Theorem \ref{theorem:T1_smooth} in the following to the particular
class of B-Spline kernels.

\begin{Definition}\label{def:spline}
  Let $\beta \in \mathbb{N}$ be given.
  The {\em normalised B-Spline} is defined by
  \begin{equation*}
    g_{\beta}\left(z\right):=\beta N_{\beta}\left(\beta z+
        \frac{\beta}{2}\right),
  \end{equation*}
  where $N_{\beta}$ denotes the cardinal B-Spline of order $\beta$.
  The cardinal B-Splines are given by $N_1(z)=1$ for $z\in(0,1)$,
  $N_1(z)=0$ elsewhere, and $N_{\beta+1}(z)=\int_{z-1}^{z}N_{\beta}(\tau)
  \d\tau$, see e.g. \cite{Chui88}.

  Furthermore, we define for $\beta \in\N$ and $N\in 2\N$ the {\em B-Spline
    kernel} by
  \begin{equation*}
    B_{\beta,N}\left(x\right)
    :=\frac{1+\eim{x}}{2\left\|g_{\beta}\right\|_{1,N}}
    \sum_{k=-\frac{N}{2}}^{\frac{N}{2}} g_{\beta}\left(\frac{k}{N}\right)
    \eip{kx}\,.
  \end{equation*}
\end{Definition}

\begin{Corollary}\label{cor:spline}
  Let $\beta \in \mathbb{N},\;\beta\ge 2$, and $N\in2\mathbb{N}$, $N\ge
  2\beta$, be given.
  Then the B-Spline kernel $B_{\beta,N}(x)$, cf. Definition \ref{def:spline},
  fulfils
  \begin{equation*}
    \left| B_{\beta,N}\left(x\right)\right| \le
    \frac{\left(2^{\beta}-1\right)\zeta\left(\beta\right) \beta^{\beta}
    }{2^{\beta-1}\pi^{\beta}-\zeta\left(\beta\right)}
    \left|Nx\right|^{-\beta}
  \end{equation*}
  for $x\in\left[-\frac{1}{2},\frac{1}{2}\right]\setminus \left\{0\right\}$
  and $B_{\beta,N}(0)=1$.
\end{Corollary}

Proof. Note that $g_{\beta} \in BV_0^{\beta-1}$.
Using $N_{\beta}'(z)=N_{\beta-1}(z)-N_{\beta-1}(z-1)$, we conclude
\begin{equation*}
  \left|g_{\beta}^{\left(\beta-1\right)}\right|_V
  =\beta^{\beta}\left|N_{\beta}^{\left(\beta-1\right)}\right|_V
  =\beta^{\beta}\left|\sum_{\tau=0}^{\beta-1} \left(-1\right)^{\tau} {\beta-1
  \choose \tau} N_1\left(\cdot-\tau\right)\right|_V
  =\left(2\beta\right)^{\beta}
\end{equation*}
and apply Theorem \ref{theorem:T1_smooth}.
\Box

\begin{figure}[ht]
  \centering
  \subfigure{\includegraphics[width=4cm]{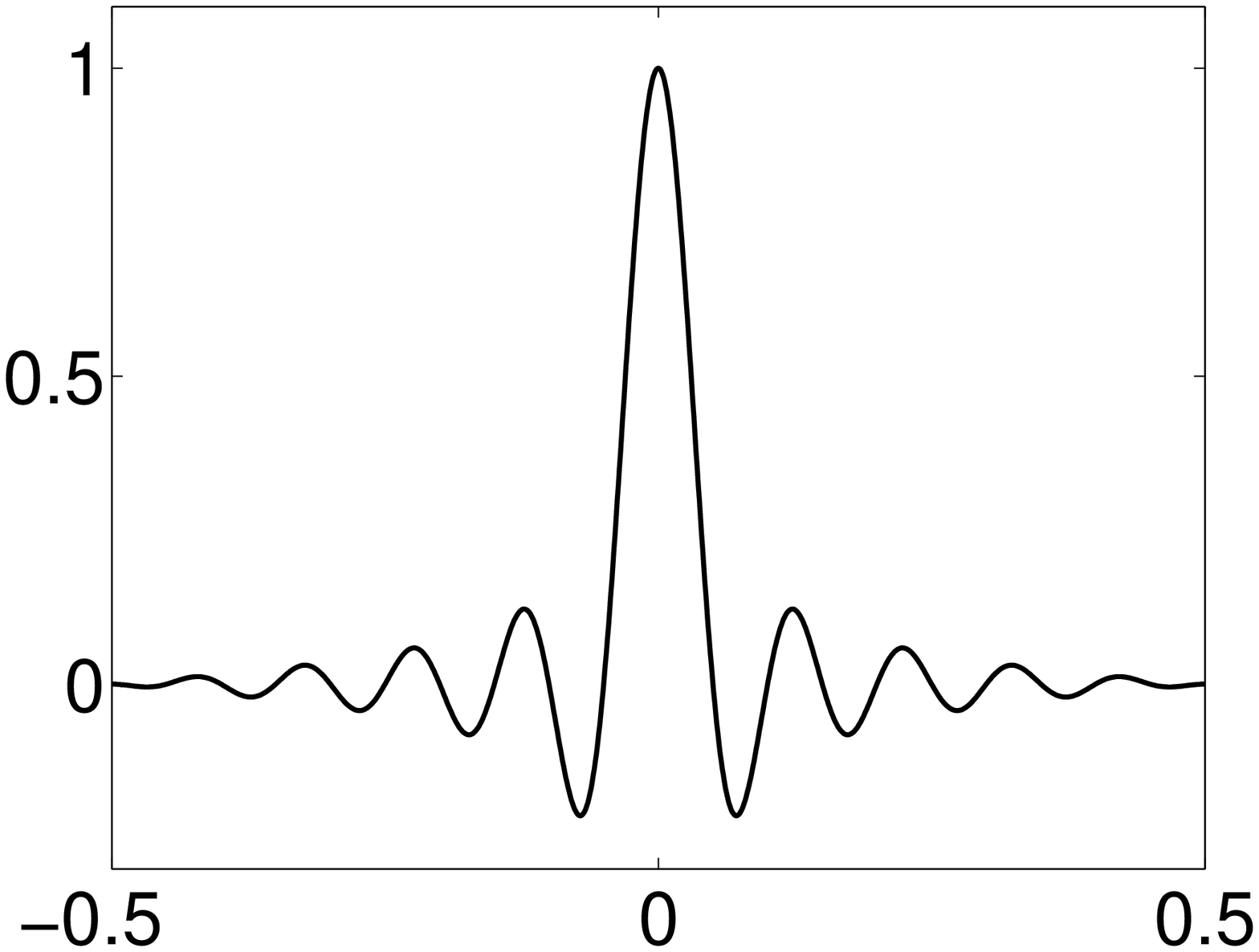}}\hfill
  \subfigure{\includegraphics[width=4cm]{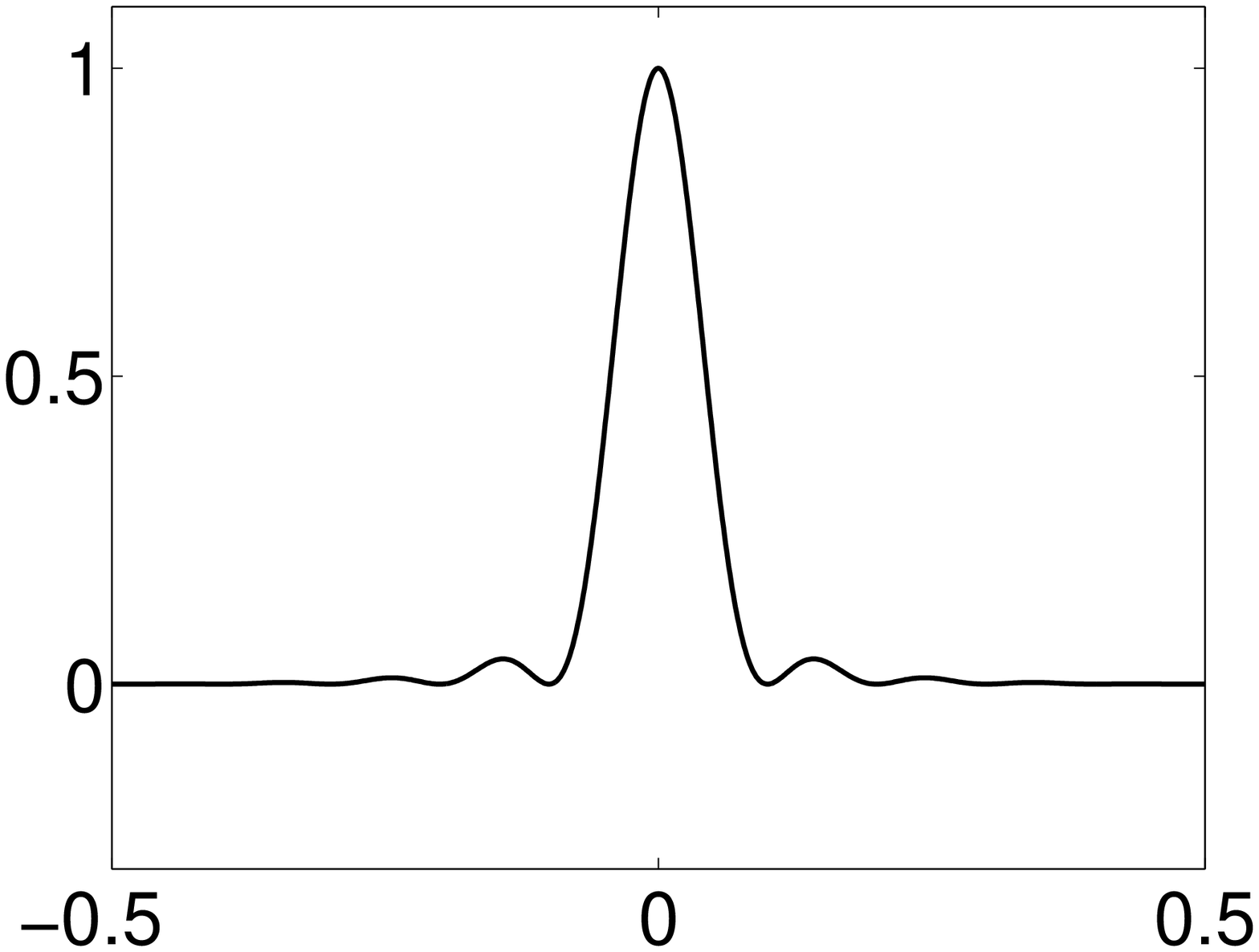}}\hfill
  \subfigure{\includegraphics[width=4cm]{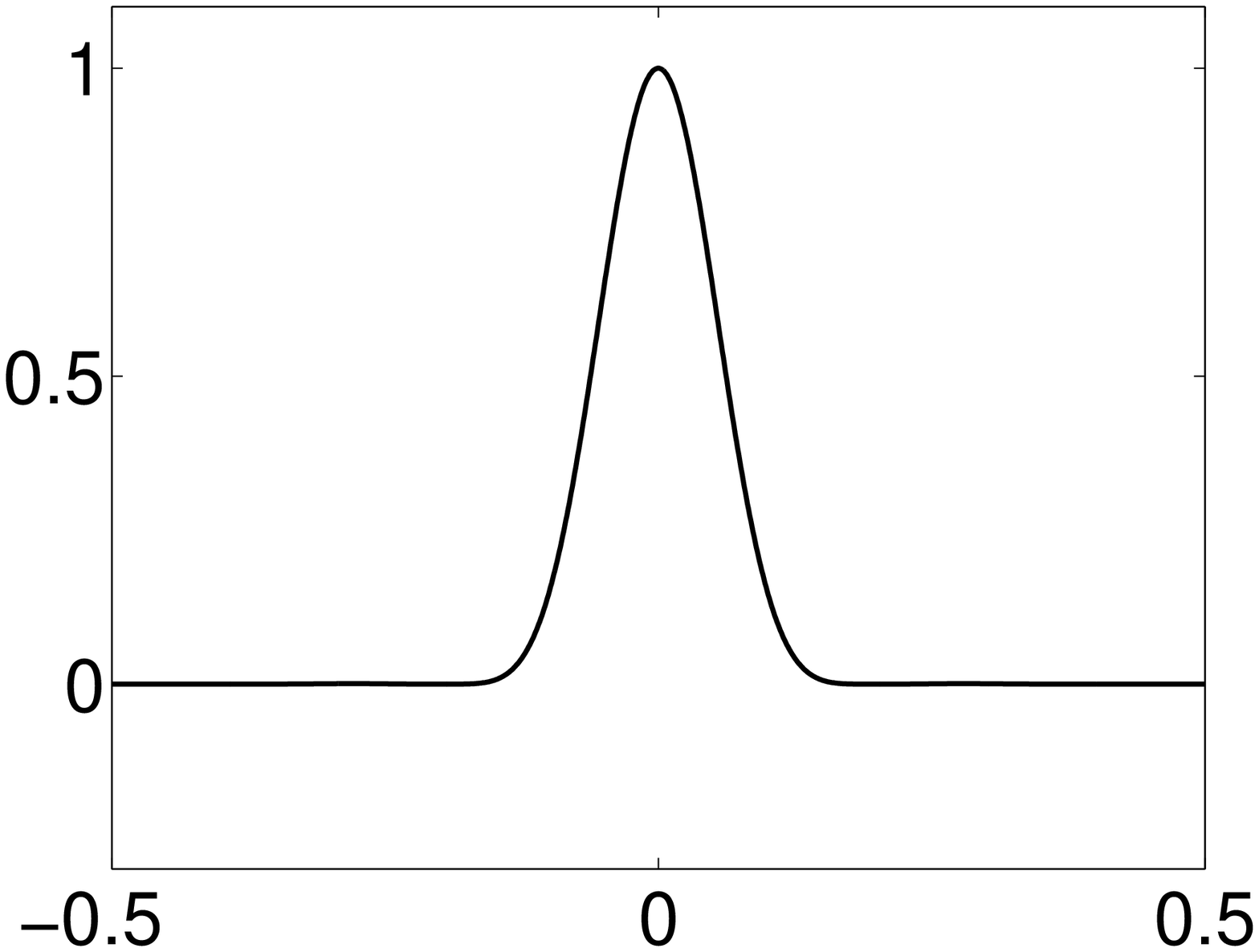}}\hfill
  \caption{From left to right: Real part of the Dirichlet kernel $B_{1,20}$,
    the Fej\'er kernel $B_{2,20}$, and the B-Spline kernel $B_{4,20}$.
    \label{Fig:polykerne}}
\end{figure}

Note, that in contrast to \cite{NaSiWa98} the order $\beta$ of the B-Spline
and the degree $N\in 2\N$ of the kernel $B_{\beta,N}$ are independent of each
other.
The special case $\beta=1$, i.e. the ``top-hat function'' $g_1(z)=1$ for
$|z|< \frac{1}{2}$ and $g_1(z)=0$ elsewhere, leads to the well known
{\em Dirichlet kernel} $B_{1,N}(x)=\frac{1}{N}\sum_{k\in\indexset}\eip{kx}$.
Analogously, $\beta=2$, i.e. the ``hat function'' $g_2(z)=2-4|z|$ for
$|z|\le \frac{1}{2}$ and $g_2(z)=0$ elsewhere, leads to the
{\em Fej\'er kernel}.
The increasing localisation of the B-Spline kernels is illustrated in Figure
\ref{Fig:polykerne}.

\begin{Remark}
  If we assume in Corollary \ref{cor:spline} furthermore, that
  $N=\beta\sigma$, $\sigma\in \N$, then the stronger estimate
  $|B_{\beta,N}(x)| \le 2\zeta(\beta)(1-2^{-\beta})
  (\frac{\beta}{\pi})^{\beta}|Nx|^{-\beta}$ holds true.
  This improvement is due to $\|g_{\beta}\|_{1,N}=N$ in Lemma
  \ref{lemma:decay_g} and follows from the
  partition of unity of the cardinal B-Spline $N_{\beta}$ and the refinement
  equation $N_{\beta}(z)=\sum_{\tau\in \Z} a_{\tau}^{(\beta,\sigma)}
  N_{\beta}(\sigma z - \tau)$ for some finitely supported coefficients
  $a_{\tau}^{(\beta,\sigma)}>0$, see e.g. \cite[pp. 8]{Chui88}.

  In particular, the Fej\'er kernel fulfils $|B_{2,N}(x)| \le |Nx|^{-2}$, $N\in
  2\N$, which also follows from the estimate $|\sin(\pi x)|\ge 2|x|$ for
  $|x|\le 1/2$ and the representation
  \begin{equation*}
    B_{2,N}\left(x\right)=\frac{2\left(1+\eim{x}\right)}{N^2}
    \left(\frac{\sin\left(\frac{N}{2} \pi x\right)}{\sin 
        \left(\pi x\right)}\right)^{2}.
  \end{equation*}

  Along the same line follows the localisation property
  $B_{1,N}(x)\le|Nx|^{-1}$
  and $B_{1,N}(0)=1$ for the Dirichlet kernel.
\end{Remark}

We complete this section by an extension of our result to the multivariate
case $d>1$.
Indeed, tensor products of the kernels constructed in Theorem
\ref{theorem:T1_smooth} yield also localised multivariate kernels as shown in
the following corollary.

\begin{Corollary}\label{cor:tensorproduct}
  Let the univariate kernel $\tilde K_N$, cf. Definition
  \ref{def:kernel}, fulfil for some $\beta\in\N$, some constant $C_{\beta}>1$,
  and $x\in\left[-\frac{1}{2}, \frac{1}{2}\right]\setminus\left\{0\right\}$
  the decay condition $|\tilde K_N(x)| \le C_{\beta} |Nx|^{-\beta}$, then its
  tensor product kernel $K_N(\zb x)=\prod_{t=0}^{d-1} \tilde K_N(x_t)$
  fulfils for $\zb x\in\left[-\frac{1}{2}, \frac{1}{2}\right]^d\setminus\left\{\zb
  0\right\}$ the estimate
  \begin{equation*}
    \left|K_N\left(\zb x\right)\right| \le \frac{C_{\beta}}
    {N^{\beta}\|\zb x\|_{\infty}^{\beta}}.
  \end{equation*}
\end{Corollary}

Proof.
The assertion follows simply from the estimate $|K_N(\zb x)| \le |\tilde
K_N(\|\zb x\|_{\infty})|$. \Box

\section{Stability of the interpolation problem}
The nonequispaced Fourier matrix $\zb A$ has full rank $M$ for $d=1$ and
$N\ge M$.
Unfortunately and due to the famous result of Mairhuber-Curtis, see e.g.
\cite[Thm. 2.3]{Wen05}, we cannot expect full rank of the matrix $\zb A$
for $N\ge M^{\frac{1}{d}}$ in the multivariate case $d>1$.

However, we prove stability results for the trigonometric interpolation
problem at $q$-separated nodes in the univariate case, cf. Theorem 
\ref{theorem:T1allgemein} and the multivariate case, cf. Theorem
\ref{theorem:Tdallgemein}.
Basically, a localised kernel $K_N$ yields a diagonal dominated kernel matrix
$\zb K_N$ and thus full rank of the nonequispaced Fourier matrix $\zb A$ for
$q$-separated nodes.
Furthermore, we prove stability results for a slightly generalised
interpolation problem at equispaced nodes and subsets of equispaced nodes in
Theorem \ref{theorem:eq}.
These results are applied to the B-Spline kernels from Section
\ref{sect:kernels} in Corollary \ref{cor:qsep} and Corollary \ref{cor:eq}.

\subsection*{The univariate setting}
The following Theorem \ref{theorem:T1allgemein} gives estimates for the
extremal eigenvalues of the matrix $\zb K_N$ under reasonable assumptions on
the kernel $K_N$.

\begin{Theorem} \label{theorem:T1allgemein}
  Let $N \in 2\mathbb{N}$ be given and let the kernel $K_N$, cf. Definition
  \ref{def:kernel}, fulfil for some $\beta> 1$ and
  $x\in\left[-\frac{1}{2},\frac{1}{2}\right]\setminus \left\{0\right\}$ 
  the localisation property
  \begin{equation*}
    \left| K_N\left(x\right)\right| \le \frac{C_{\beta}}{ N^{\beta}
      |x|^{\beta}}\,.
  \end{equation*}
  Furthermore, let a sampling set ${\cal X}$ contain arbitrary nodes with
  separation distance $q>0$.
  Then, the extremal eigenvalues of the matrix $\zb K_N$ are bounded by
  \begin{equation*}
    1 -  \frac{2\, \zeta\left(\beta\right) C_{\beta}}{N^{\beta} q^{\beta}}
    \le \lambda \le 1 \le \Lambda \le
    1 + \frac{2\, \zeta\left(\beta\right) C_{\beta}}{N^{\beta}  q^{\beta}}.
  \end{equation*}
\end{Theorem}

Proof.
As usual, let $M$ denote the number of nodes in $\cal X$.
Due to $K_N\left(0\right)=1$, cf. Definition \ref{def:kernel}, we obtain
$\trace(\zb K_N):=\sum_{j=0}^{M-1} K_N(0)=M$.
Since the trace is invariant under similarity transforms, all eigenvalues
sum up to $M$ and thus, the inequality $\lambda \le 1 \le \Lambda$ is
fulfilled.

Now, let $\lambda_{\star}$ be an arbitrary eigenvalue of $\zb K_N$, then for
some index $j\in\left\{0,\hdots,M-1\right\}$ Gershgorin's circle theorem
yields
\begin{equation*}
  \left|\lambda_{\star}-1\right|
  \le \sum\limits_{l=0; l\ne j}^{M-1}
  \left|K_N\left(x_j-x_l\right)\right|.
\end{equation*}
By using that the separation distance of the sampling set is $q$ and by the
localisation of the kernel $K_N$, we obtain
\begin{eqnarray*}
  \left|\lambda_{\star}-1\right|
  &\le & \frac{C_{\beta}}{N^{\beta} }\sum\limits_{l=0; l\ne j}^{M-1}
  \frac{1}{\left| x_j-x_l\right|^\beta}
  \le  \frac{2 \, C_{\beta}}{N^{\beta} q^{\beta}}
  \sum\limits_{l=1}^{\left\lfloor M/2\right\rfloor} l^{-\beta}
    <   \frac{2\,\zeta(\beta) C_{\beta}}{N^{\beta} q^{\beta}}.
\end{eqnarray*}
\Boxgl

Indeed, the kernels constructed in Theorem \ref{theorem:T1_smooth} yield well
conditioned matrices $\zb K_N$.
We also note that the decay of the Dirichlet kernel only allows for a weaker
result.
\begin{Remark}
  Let a sampling set ${\cal X}\in\T$ with separation distance
  $q>0$ be given.
  Then the application of Theorem \ref{theorem:T1allgemein}, where the last
  step of its proof is replaced by the estimate $\sum_{l=1}^{\left\lfloor
    M/2\right\rfloor} l^{-1} \le 1+\ln\frac{1}{2q}$, yields:
  The matrix $(B_{1,N}(x_j-x_l)_{j,l=0,\hdots,M-1}=\frac{1}{N}
  \zb A\zb A^{\adj}$ is nonsingular for $N>\left(1+\left|\log\left(2q\right)
  \right|\right)q^{-1}$.
  The logarithmic term in this condition is clearly suboptimal.
\end{Remark}

As an immediate consequence of Theorem \ref{theorem:T1allgemein} we state a
stability result for an equispaced grid disturbed by jitter.
\begin{Corollary}\label{cor:jitter}
  Let the assumptions of Theorem \ref{theorem:T1allgemein} hold true.
  Furthermore, let the sampling nodes be of the form
  $x_j=-\frac{1}{2}+\frac{j-\varepsilon_j}{M},\; j=0,\hdots,M-1,$ where
  $0\le\varepsilon_j \le \varepsilon < 1$.
  Then the eigenvalues of the matrix $\zb K_N$ are bounded by
  \begin{equation*}
    1-\frac{2\, \zeta\left(\beta\right) C_{\beta} M^{\beta}}{N^{\beta} \left(1-
        \varepsilon \right)^{\beta}} \le \lambda \le 1 \le \Lambda \le
    1+\frac{2\, \zeta\left(\beta\right) C_{\beta} M^{\beta}}{N^{\beta} \left(1-
        \varepsilon \right)^{\beta}}\,.
  \end{equation*}
\end{Corollary}

Proof. Since the separation distance is bounded by $q\ge M^{-1}
(1-\varepsilon)$ the result follows by Theorem \ref{theorem:T1allgemein}. \Box

\subsection*{The multivariate setting}
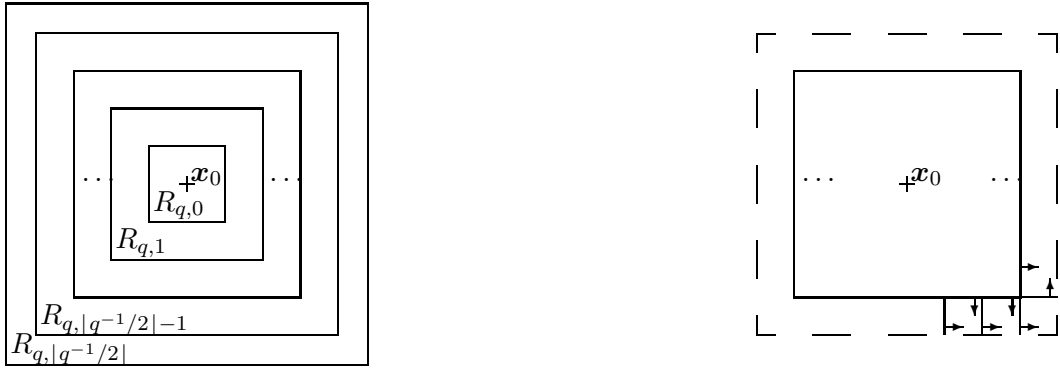
\begin{figure}[ht!]
  \centering
  \subfigure{\input{images/figureTd}} \hfill
  \subfigure{\input{images/figureTd_refined}}
  \caption{Partitioning of the torus $\mathbb{T}^2$ into the rings
    $R_{q,m},\;m=0,\hdots\lfloor q^{-1}/2 \rfloor$ (left).
    Further subdivision into shifted and rotated versions of the cube
    $[0,q)^d$, where arrows indicate the ``ownership'' of the faces to a
    particular cube (right).\label{fig:Td}}
\end{figure}

First, we borrow a packing argument on the sphere from \cite{NaSiWa98} and
refine it in Lemma \ref{lemma:Tda} for the present setting, i.e., we show how
many $q$-separated nodes can be placed in a certain distance to a reference
node, see also Figure \ref{fig:Td}.

\begin{Definition}\label{def:part}
  For $d\in\N$ and a separation distance $q$, $0<q \le \frac{1}{2}$, we define
  the partitioning
  \begin{equation*}
    R_{q,m}:=
    \left\{\zb x \in \T^d:m q \le \dist\left(\zb x,\zb
    0\right)<\left(m+1\right)q\right\}
  \end{equation*}
  for $m=0,\hdots,\left\lfloor q^{-1}/2 \right\rfloor-1$ and
  \begin{equation*}
    R_{q,\left\lfloor q^{-1}/2 \right\rfloor}:=
    \left\{\zb x \in \T^d: \left\lfloor q^{-1}/2 \right\rfloor q \le
    \dist\left(\zb x,\zb 0\right)\le 1/2\right\}.
  \end{equation*}
  Its restriction to the sampling set ${\cal X}$ will be denoted by $R_{{\cal
      X},q,m}:=R_{q,m}\cap {\cal X}$.
\end{Definition}

\begin{Lemma} \label{lemma:Tda}
   Let $d\in \N$ and an $q$-separated sampling set ${\cal X}$ with $0<q \le
   \frac{1}{2}$ be given.
   Then, each of the sets $R_{{\cal X},q,m}$ has bounded cardinality
   \begin{equation*}
     \left|R_{{\cal X},q,m}\right|\le 2^d\left(2^d-1\right) m^{d-1},
     \qquad m=1,\hdots,\left\lfloor q^{-1}/2 \right\rfloor.
   \end{equation*}
\end{Lemma}
Proof. 
We use a packing argument for the partition $\{R_{q,m},\,m=0,\hdots,\lfloor
q^{-1}/2 \rfloor\}$ of the torus $\T^d$.
Each ring $R_{q,m}$ is subdivided into shifted and rotated versions of the
cube $[0,q)^d$, cf. Figure \ref{fig:Td} (right).
This is done such that each point in $R_{q,m}$ is contained in at least one
of these boxes and the boxes share no interior points with each other.
Every box contains at most one node of the sampling set and hence, the
estimate
\begin{equation*}
  \left|R_{{\cal X},q,m}\right| \le \frac{1}{q^d}
  \int\limits_{R_{q,m} } \d \zb x \le
  2^d\left(\left(m+1\right)^d-m^d\right)
  =2^d\sum_{t=1}^d {d\choose t} m^{d-t}
  \le 2^d m^{d-1} \sum_{t=1}^d {d\choose t}.
\end{equation*}
is valid.\Box

Using localised kernels in conjunction with a separated sampling set and Lemma
\ref{lemma:Tda}, we state the following theorem on the stability of the
interpolation problem.

\begin{Theorem} \label{theorem:Tdallgemein}
  Let $d\in\N$, $N \in 2\mathbb{N}$ be given and let the kernel $K_N$,
  cf. Definition \ref{def:kernel}, fulfil for some $\beta> d$ and
  $\zb x\in\left[-\frac{1}{2},\frac{1}{2}\right]^d\setminus\left\{\zb
  0\right\}$ the localisation property
  \begin{equation*}
    \left|K_N\left(\zb x\right)\right| \le \frac{C_{\beta}}
    {N^{\beta}\|\zb x\|_{\infty}^{\beta}}\,.
  \end{equation*}
  Furthermore, let a sampling set ${\cal X}$ contain arbitrary nodes with
  separation distance $0<q\le \frac{1}{2}$.
  Then, the extremal eigenvalues of the matrix $\zb K_N$ are bounded by
  \begin{equation*}
    1 - \frac{2^d\left(2^d-1\right) \zeta\left(\beta-d+1\right)
    C_{\beta}}{N^{\beta} q^{\beta}} \le \lambda \le 1 \le \Lambda \le
    1 + \frac{2^d\left(2^d-1\right) \zeta\left(\beta-d+1\right)
    C_{\beta}}{N^{\beta} q^{\beta}}.
  \end{equation*}
\end{Theorem}

Proof. 
Let $\lambda_{\star}$ be an arbitrary eigenvalue of $\zb K_N$.
Without loss of generality, let the diagonal element of the matrix $\zb K_N$
used in Gershgorin's circle theorem correspond to $\zb x_0=\zb 0$.
Then we conclude by $K_N(\zb 0)=1$, cf. Definition \ref{def:kernel}, that
\begin{eqnarray*}
  \left|\lambda_{\star}-1\right| &\le& \sum\limits_{l=1}^{M-1} \left|
    K_N\left(\zb 0-\zb x_l\right)\right|.
\end{eqnarray*}
Using the partition from Definition \ref{def:part}, Lemma \ref{lemma:Tda},
and the localisation of the kernel $K_N$, we get
 \begin{eqnarray*}
 \left|\lambda_{\star}-1\right|
     &\le&\sum\limits_{m=1}^{\left\lfloor q^{-1}/2 \right\rfloor}
     \sum\limits_{\zb x_l \in R_{{\cal X},q,m}} 
     \left|K_N\left(-\zb x_l\right)\right|\\
     &\le& \frac{2^d\left(2^d-1\right) C_{\beta}}{N^{\beta}}
     \sum\limits_{m=1}^{\left\lfloor q^{-1}/2 \right\rfloor} m^{d-1}
     \max\limits_{\zb x \in R_{q,m}} \|\zb x\|_{\infty}^{-\beta}\\ 
     &\le& \frac{2^d\left(2^d-1\right) \zeta\left(\beta-d+1\right)
     C_{\beta}}{N^{\beta} q^{\beta}}.
   \end{eqnarray*}
\Boxgl

Particularly, this result includes Theorem \ref{theorem:T1allgemein} if we set
$d=1$. 

\begin{Corollary}  \label{cor:qsep}
  Let the dimension $d\in\N$, an arbitrary sampling set ${\cal X}\in \T^d$
  with separation distance $0<q\le\frac{1}{2}$, and a polynomial degree $N\in
  2\N$, $N>2dq^{-1}$, be given.
  Then the nonequispaced Fourier matrix $\zb A$ has full rank.
  Moreover, the eigenvalues of the kernel matrix $\zb K_N=\zb A\zb {\hat W}
  \zb A^{\adj}$ obtained from the B-Spline kernel of order $\beta=d+1$ are
  bounded by
    \begin{equation*}
      0<1- \left(\frac{2d}{Nq}\right)^{d+1}
      \le \lambda \le 1 \le \Lambda \le
      1+ \left(\frac{2d}{Nq}\right)^{d+1}.
    \end{equation*}
\end{Corollary}

Proof. Note first that $N\ge 2\beta$.
We apply Theorem \ref{theorem:Tdallgemein} where we use the estimates
for $C_{\beta}$ given in Corollary \ref{cor:spline} and simplify the involved
constant.
Hence, the full rank of $\zb A$ follows.
\Box

Thus, we have shown that the optimal trigonometric interpolation problem at
$q$-separated nodes in $d$ dimensions obeys a uniformly bounded condition
number for a polynomial degree $N > 2d q^{-1}$ and appropriate damping
factors.
The dependence on $q^{-1}$ is optimal as the subsequent analysis of the
equispaced case shows.
However, the constant $2 d$ is not optimal for high spatial dimensions.
As pointed out in \cite{BaGr03} for the related approximation problem
\eqref{eq:ap}, it is an open problem to improve on this.

\medskip
In summary, Lemma \ref{lemma:cg_error} and Corollary \ref{cor:qsep} assure in
our situation a prescribed reduction of the error $\|\zb {\hat e}_l\|_{\zb
{\hat W}^{-1}}$ in a constant number of iterations.
Hence, if we assume an additional uniformity condition $q=c M^{-\frac{1}{d}}$
for the sampling set  ${\cal X}$, the total arithmetical complexity of
Algorithm \ref{algo:CGNE} for solving \eqref{eq:optInterp} up to a prescribed
error is bounded by ${\cal O}(M\log M)$.

\begin{Remark}
Other frequently applied kernels also possess specific localisation properties
and thus, yield stable interpolation at $q$-separated nodes by means of
Theorem \ref{theorem:Tdallgemein} as follows:

We define for $\beta\in2\N$, $\sigma\in \N$, and $N=\beta(\sigma-1)+2$ the
Jackson kernel by
\begin{equation*}
  J_{\beta,N}\left(x\right):=\frac{1+\eim{x}}{2\sigma^{\beta}}
  \left(\frac{\sin\left(\sigma\pi x\right)}{\sin 
      \left(\pi x\right)}\right)^{\beta}.
\end{equation*}
Being a normalised power of the Fej\'er kernel $B_{2,2\sigma}$, the
coefficients $\hat w_k$ of the Jackson kernel can be obtained by an iterated
discrete convolution of the coefficients of the Fej\'er kernel, see
\cite{Alekseev} for details.
In contrast, the B-Spline kernel relies on a continuous convolution.
The Jackson kernel is localised as $|J_{\beta,N}(x)| \le (\frac{\beta}{2} 
)^{\beta} |Nx|^{-\beta}$ for $x\in\left[-\frac{1}{2},\frac{1}{2}\right]
\setminus\left\{0\right\}$ and fulfils $J_{\beta,N}\left(0\right)=1$.
Hence, the tensor product Jackson kernel of appropriate order yields for
$N\in 2\N$, $N>2.1 d q^{-1}$, the nonsingular kernel matrix 
\begin{equation*}
  \left(J_{2\lceil\frac{d+1}{2}\rceil,N}\left(\zb x_j-\zb x_l\right)\right)_{
    j,l=0,\hdots,M-1}.
\end{equation*}

Secondly, it is well known that the weight $1+(2\pi k)^{2\alpha}$ is
associated to the squared Sobolev norm $\|f\|_2^2+\|f^{(\alpha)}\|_2^2$.
For $\beta\in\N$ and $\alpha,\gamma>0$, a regularised and slightly generalised
weight is given by
\begin{equation*}
  g_{\alpha,\beta,\gamma}\left(z\right):=c_{\alpha,\beta,\gamma}
  \frac{\left(\frac{1}{4}-z^2\right)^{\beta}}{\gamma+|z|^{2\alpha}}
\end{equation*}
for $|z|\le\frac{1}{2}$ and $g_{\alpha,\beta,\gamma}(z)=0$
elsewhere, where the constant $c_{\alpha,\beta,\gamma}$ is chosen such
that $\|g_{\alpha,\beta,\gamma}\|_{L_1}=1$.
Here, the denominator generalises the weight $1+(2\pi k)^{2\alpha}$ and the
nominator ensures $g_{\alpha,\beta,\gamma} \in BV^{\beta-1}_0$.
We define for $N\in2\N$ the Sobolev kernel by
\begin{equation*}
  S_{\alpha,\beta,\gamma,N}\left(x\right)
  :=\frac{1+\eim{x}}{2\left\|g_{\alpha,\beta,\gamma}\right\|_{1,N}}
  \sum_{k=-\frac{N}{2}}^{\frac{N}{2}}
  g_{\alpha,\beta,\gamma}\left(\frac{k}{N}\right)
  \eip{kx}\,.
\end{equation*}
The kernel is localised as $|S_{\alpha,\beta,\gamma,N}(x)| \le \tilde
c_{\alpha,\beta,\gamma} |Nx|^{-\beta}$ for
$x\in\left[-\frac{1}{2},\frac{1}{2}\right] \setminus\left\{0\right\}$
and some constant $\tilde c_{\alpha,\beta,\gamma}>0$ and fulfils
$S_{\alpha,\beta,\gamma,N}\left(0\right)=1$.
\end{Remark}

\subsection*{Results for equispaced nodes} \label{sect:eq}
In the case of equispaced nodes we employ the fact that the matrix $\zb K_N$
is circulant.
We present a slightly generalised result in the following Theorem
\ref{theorem:eq}.

\begin{Definition}  \label{def:l1kernel}
  We define for $d,n\in\N$ and weights $\hat w_{\zb k} \in \R,\; \zb k \in
  \Z^d$ with $\sum_{\zb k\in \Z^d}|\hat w_{\zb k}| < \infty$, the kernel
  \begin{equation*}
    K \left(\zb x\right):=\sum\limits_{\zb k \in \Z^d} \hat w_{\zb k}
    \eip{\zb k \zb x}
  \end{equation*}
  and by evaluating at the equispaced sampling nodes $\zb j = (j_0,\hdots, 
  j_{d-1})^{\top} \in I_n$, the matrix
  \begin{equation*}
    \zb K: = \left(K \left(\frac{\zb j-\zb l}{n}\right) \right)_{\zb j,\zb l
      \in I_n} \in \C^{n^d\times n^d}.
  \end{equation*}
\end{Definition}

\begin{Theorem}\label{theorem:eq}
  The matrix $\zb K$ in Definition \ref{def:l1kernel} possesses the following
  properties.
  Its eigenvalues are given by
  \begin{equation*}
    \lambda_{\zb s}\left(\zb K\right)
    = n^d \sum_{\zb r\in \Z^d} \hat w_{\zb s+n\zb r}
  \end{equation*}
  for $\zb s\in I_n$.
  For tensor product weights $\hat w_{\zb k}=\prod_{t=0}^{d-1} \hat w_{k_t}$,
  $\sum_{k\in\Z} |\hat w_k|<\infty$, this simplifies to
  \begin{equation*}
    \lambda_{\zb s}\left(\zb K\right)
    = n^d \prod_{t=0}^{d-1} \sum_{r_t\in \Z} \hat w_{s_t+n r_t}
  \end{equation*}
  for $\zb s\in I_n$.
  Moreover, the extremal eigenvalues of
  \begin{equation*}
    \zb K_{\Gamma}: = \left(K\left(\frac{\zb j - \zb l}{n}\right) \right)_{
      \zb j, \zb l\in \Gamma}
  \end{equation*}
  are bounded by the extremal eigenvalues of $\zb K$ for any $\Gamma\subset
  I_n$. 
\end{Theorem}
Proof.
  The matrix $\zb K$ is multilevel circulant and thus diagonalised by the
  Fourier matrix $\zb F_n=(\eip{\zb k \zb j/n})_{\zb j,\zb k\in I_n}$.
  We calculate
  \begin{eqnarray*}
    \left(\zb F_n^{\adj} \zb K \zb F_n\right)_{\zb s,\zb t}
    &=&\sum_{\zb j,\zb l\in I_n} \eim{\zb s \zb j/n} K\left(\frac{\zb j-
    \zb l}{n}\right) \eip{\zb t \zb l/n}\\
    &=&\sum_{\zb k\in\Z^d} \hat w_{\zb k}
       \sum_{\zb j\in I_n} \eim{\zb j\left(\zb s-\zb k\right)/n}
       \sum_{\zb l\in I_n} \eip{\zb l\left(\zb t-\zb k\right)/n}
  \end{eqnarray*}
  for $\zb s,\zb t\in I_n$ and use
  \begin{equation*}
    \sum_{\zb j\in I_n} \eip{\zb j\left(\zb s-\zb k\right)/n}=
    \begin{cases}
      n^d & \text{if} \quad \frac{\zb s-\zb k}{n}\in\Z^d,\\
      0                & \text{otherwise}.
    \end{cases}
  \end{equation*}
  See also \cite[Cor. 3.10, Thm. 3.11]{NaSiWa98} for the univariate case.

  The second assertion is due to the Kronecker product structure of the
  matrix $\zb K$ in the case of tensor product kernels.
  The last assertion follows from the fact that removing a node is nothing
  else than removing its corresponding row and column in $\zb K$ and from the
  interlacing property for eigenvalues, see \cite[pp. 185]{hojo}.
\Box

Now, let $d,n\in\N$, $N\in 2\N$, and the equispaced sampling set ${\cal X}=
\frac{1}{n} I_n\subset\T^d$ be given.
A simple consequence of Theorem \ref{theorem:eq} is the fact that the kernel
matrix $\zb K_N=\frac{1}{N}\zb A\zb A^{\adj}$ is singular whenever
$N<q^{-1}=n$.
Hence, the condition $N>2dq^{-1}$ for the full rank of $\zb A$ is optimal with
respect to $q$.
Nevertheless, we also apply Theorem \ref{theorem:eq} to obtain positive
results in the equispaced setting for the Dirichlet and the Fej\'er kernel.
\begin{Corollary}  \label{cor:eq}
  Let the dimension $d\in\N$, a sampling set $\cal X$ with $M\ge 2$ equispaced
  nodes, and a polynomial degree $N\in 2\mathbb{N}$ with
  $N>q^{-1}=M^{\frac{1}{d}}$ be given.
  Then, for the Dirichlet kernel $B_{1,N}$ and its tensor product version for
  $d>1$, the extremal eigenvalues of the kernel matrix $\zb K_N$ obey
  \begin{equation*}
    \left(\frac{\left\lfloor Nq \right\rfloor}{Nq}\right)^d
    = \lambda \le 1 \le \Lambda =
    \left(\frac{\left\lceil{Nq}\right\rceil}{Nq}\right)^d.
  \end{equation*}
  Furthermore, the Fej\'er kernel $B_{2,N}$ yields
  \begin{equation*}
    \left(1- \frac{1}{ N^2q^{2}}\right)^d
    \le \lambda \le 1 \le \Lambda \le
    \left(1+ \frac{1}{ N^2q^{2}}\right)^d
  \end{equation*}
  where equality holds for the outmost inequalities if
  $N=\left(2\sigma+1\right)q^{-1},\,\sigma\in\mathbb{N}$.
  In particular, the kernel matrix $K_N$ is nonsingular.
\end{Corollary}

Proof. Throughout this proof, let $n=q^{-1}$ and the damping factors be
extended by $\hat w_{\zb k}=0$ for $\zb k \notin \indexset$.
We apply for $d=1$ the first statement of Theorem \ref{theorem:eq} to the
weights $\hat w_k= \frac{1}{N}$, $k \in I_N$, of the Dirichlet kernel for
$N\in 2\N$ and to the weights $\hat w_k= \frac{2}{N}(1-\frac{|2k+1|}{N})$,
$k \in I_N$, of the Fej\'er kernel, if $N=\left(2\sigma+1\right)q^{-1},\,
\sigma\in \mathbb{N}$, respectively.

The assertion is little more delicate for the univariate Fej\'er kernel and
$N\neq\left(2\sigma+1\right)q^{-1}$.
We use the representation
\begin{equation*}
  \left|B_{2,N}\left(x\right)\right| = \frac{4}{N^2}
  \sum\limits_{r=0}^{\frac{N}{2}-1} \sum\limits_{k=-r}^{r} \eip{kx}.
\end{equation*}
Now, let $\lambda_{\star}$ be an arbitrary eigenvalue of the kernel matrix
$\zb K_N$, then Gershgorin's circle theorem yields
\begin{equation*}
  \left|\lambda_{\star}-1\right|
  \le \sum_{l=1}^{n-1}\left|B_{2,N}\left(\frac{l}{n}\right)\right|
   = \frac{4n}{N^2} \sum\limits_{r=0}^{\frac{N}{2}-1}
   \left(2\left\lfloor \frac{r}{n} \right\rfloor+1\right) - 1.
\end{equation*}
Since for $Q:=\left\lfloor\frac{N-2}{2n}\right\rfloor$ and $R:=\frac{N}{2}-1-nQ$ the identity
\begin{equation*}
  \sum\limits_{r=0}^{\frac{N}{2}-1}\left\lfloor \frac{r}{n} \right\rfloor
  =\sum\limits_{s=0}^{Q-1}\;\sum\limits_{r=sn}^{\left(s+1\right)n-1} s + \sum_{r=nQ}^{nQ+R} Q
  =\frac{\left(N-n\right)^2 - \left(2\left(R+1\right)-n\right)^2}{8n}
\end{equation*}
holds, we proceed
\begin{eqnarray*}
  \frac{4n}{N^2} \sum\limits_{r=0}^{\frac{N}{2}-1}
  \left(2\left\lfloor \frac{r}{n} \right\rfloor+1\right) - 1
  &= &\frac{4n}{N^2} \left(2 \frac{\left(n-N\right)^2 - \left(2\left(R+1\right)-n\right)^2}{8n} + \frac{N}{2}\right) - 1\\
  &= &\frac{n^2}{N^2}-\left(\frac{2\left(R+1\right)-n}{N}\right)^2\\
  &\le & \frac{n^2}{N^2}.
\end{eqnarray*}

The case $d>1$ is due to the second statement in Theorem.
\Box


\section{Numerical results}
In this section, we exemplify our findings on the stability of the optimal
interpolation problem \eqref{eq:optInterp} and its iterative solution by
Algorithm \ref{algo:CGNE}.

\begin{figure}[ht!]
  \centering
  \subfigure{\includegraphics[width=6cm]{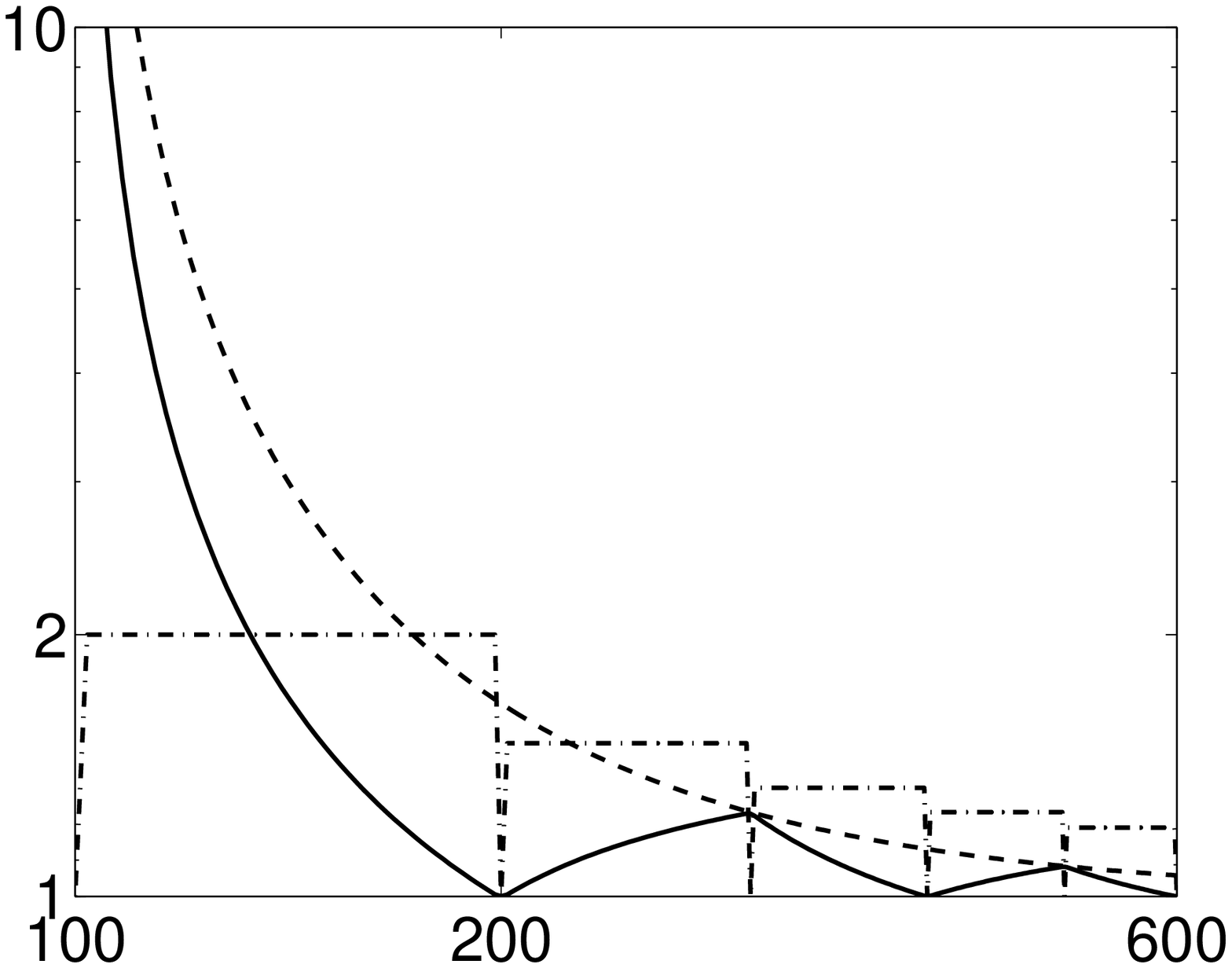}}\hfill
  \subfigure{\includegraphics[width=6cm]{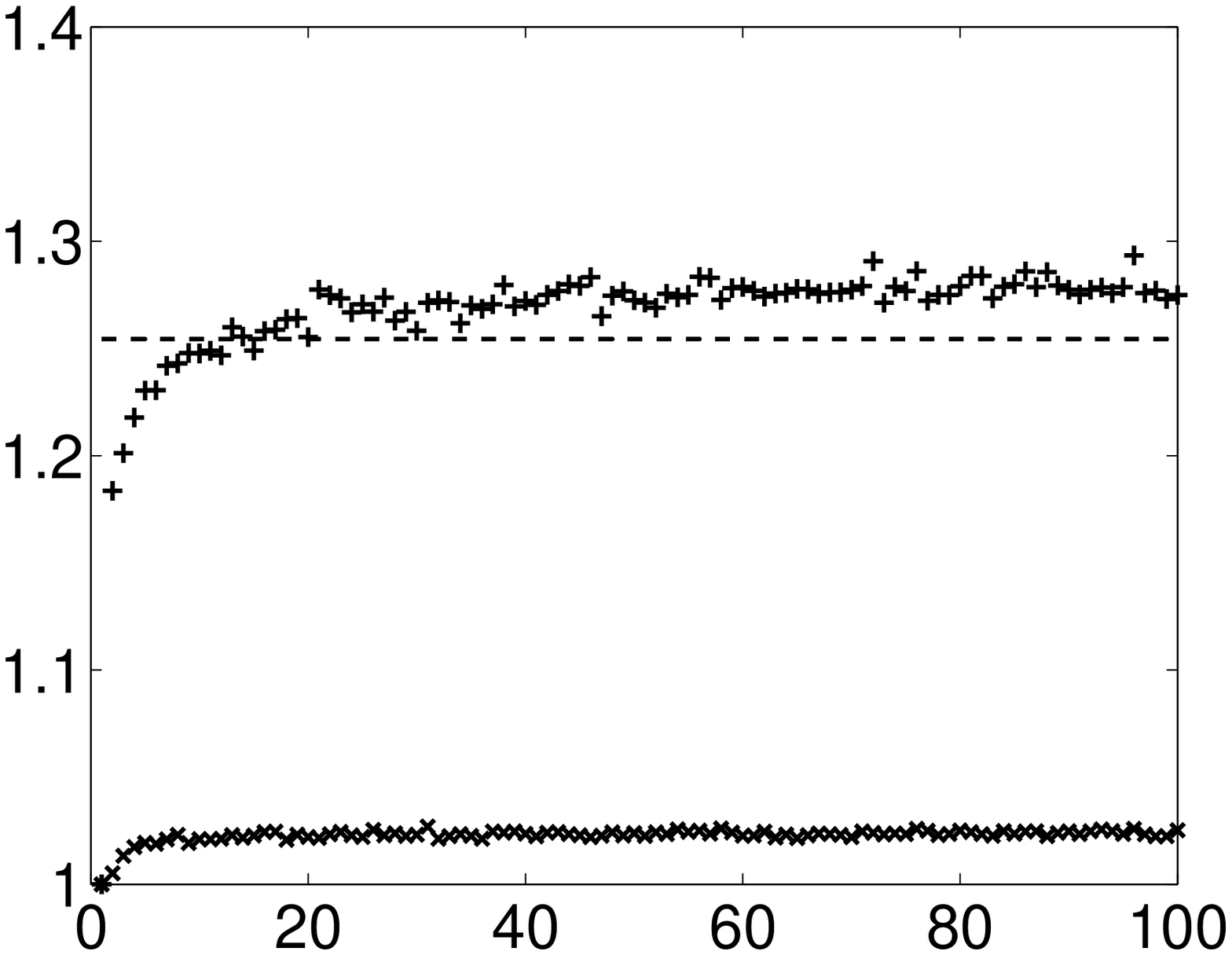}}
  \caption{Condition number of the kernel matrix $\zb K_N$.
    Left: Condition number with respect to polynomial degree
    $N=100,\hdots,600$, no weights, i.e., Dirichlet kernel (dash-dot); weight
    function $g_2$, i.e., Fej\'er kernel (solid), and the estimate of Corollary
    \ref{cor:eq} (dashed); here, the number of equispaced nodes is $M=100$.
    Right: Condition number with respect to the number of nodes
    $M=1,\hdots,100$, the nodes are equispaced perturbed by
    $\varepsilon_{\text{rel.}}=0.1$ jitter error, the polynomial degree is 
    $N=6M$;
    no weights, i.e., Dirichlet kernel ($\,+\,$);
    weight function $g_2$, i.e., Fej\'er kernel (${}\times{}$), and its estimate
    by Corollary \ref{cor:jitter} (dashed).\label{Fig:cond_eq_jitter}}
\end{figure}

The estimates for the condition number of the kernel matrix $\zb K_N$ for
equispaced nodes, cf. Corollary \ref{cor:eq}, are shown in Figure
\ref{Fig:cond_eq_jitter} (left).
For $Nq\in \N$ and the Dirichlet kernel $B_{1,N}$ the matrix $\zb K_N$ is
just the identity.
However, using the better localised Fej\'er kernel $B_{2,N}$ improves the
condition number already for $N>\sqrt{3}q^{-1}$ when $Nq\not\in \N$.

We present the effect on the stability of the interpolation problem when the
equispaced nodes are perturbed by jitter error, cf. Corollary
\ref{cor:jitter}, in Figure \ref{Fig:cond_eq_jitter} (right).
We choose different sampling sets of size $M=1,\hdots,100$ with equispaced
nodes disturbed by $10\%$ jitter error and evaluate the maximum condition
number over $100$ reruns for the Dirichlet kernel $B_{1,6M}$ and the Fej\'er
kernel $B_{2,6M}$, respectively.
The Fej\'er kernel produces a lower condition number which is also validated by
the shown upper bound. These results confirm the theoretical results
of Corollary \ref{cor:eq} and Corollary \ref{cor:jitter}.

\begin{figure}[ht]
  \centering
  \subfigure{\includegraphics[width=6cm]{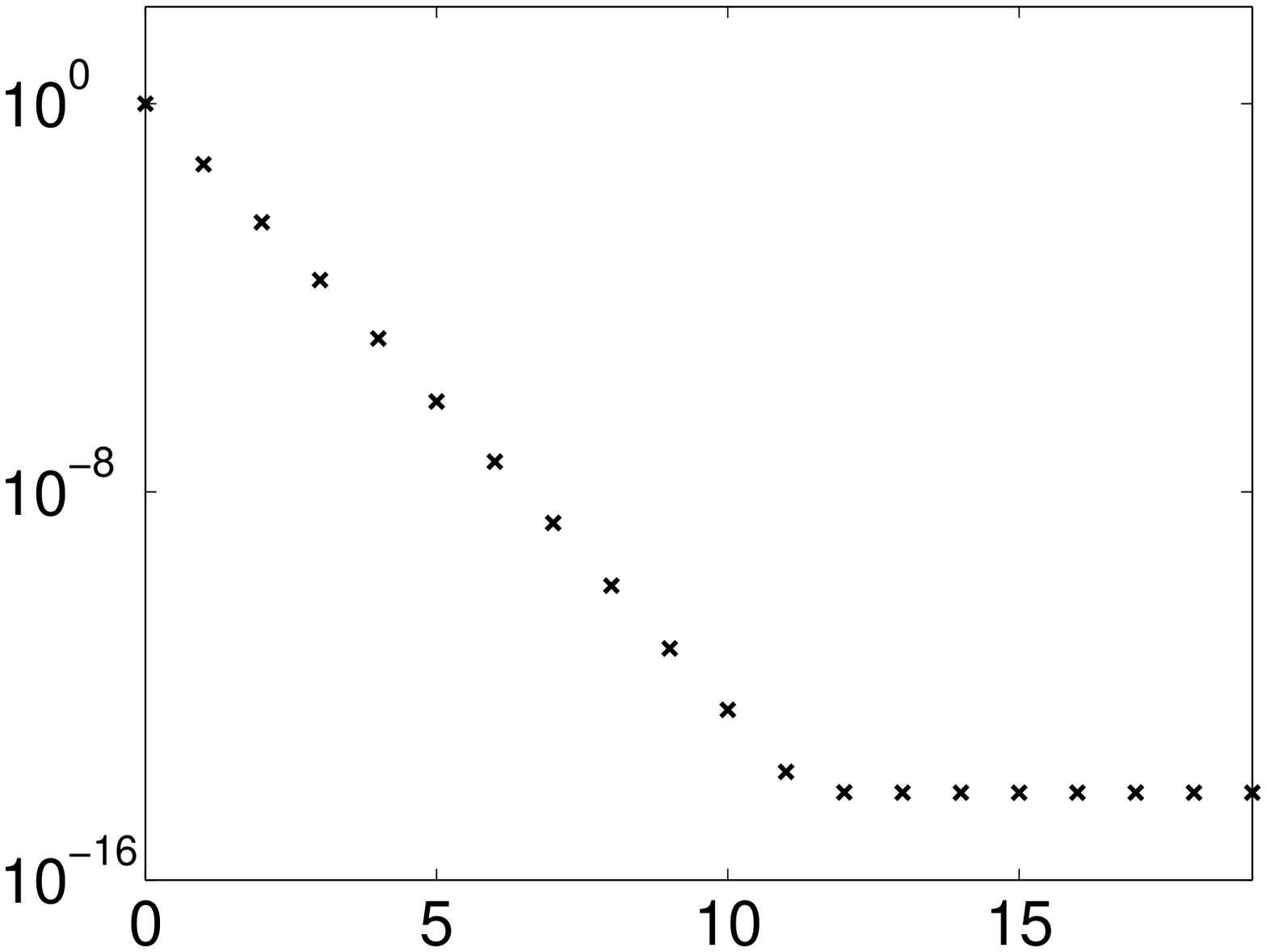}}\hfill
  \subfigure{\includegraphics[width=6cm]{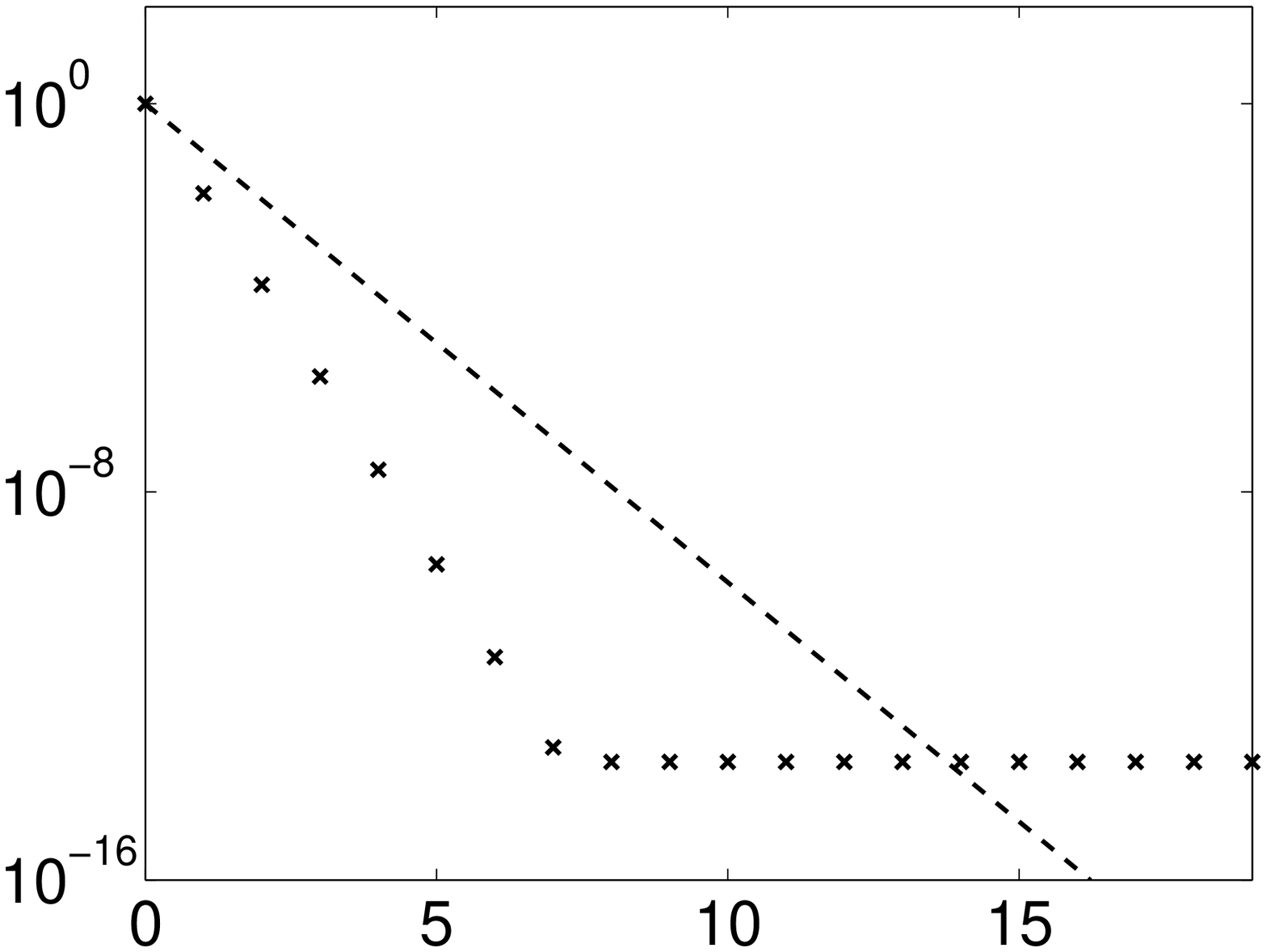}} \\
  \subfigure{\includegraphics[width=6cm]{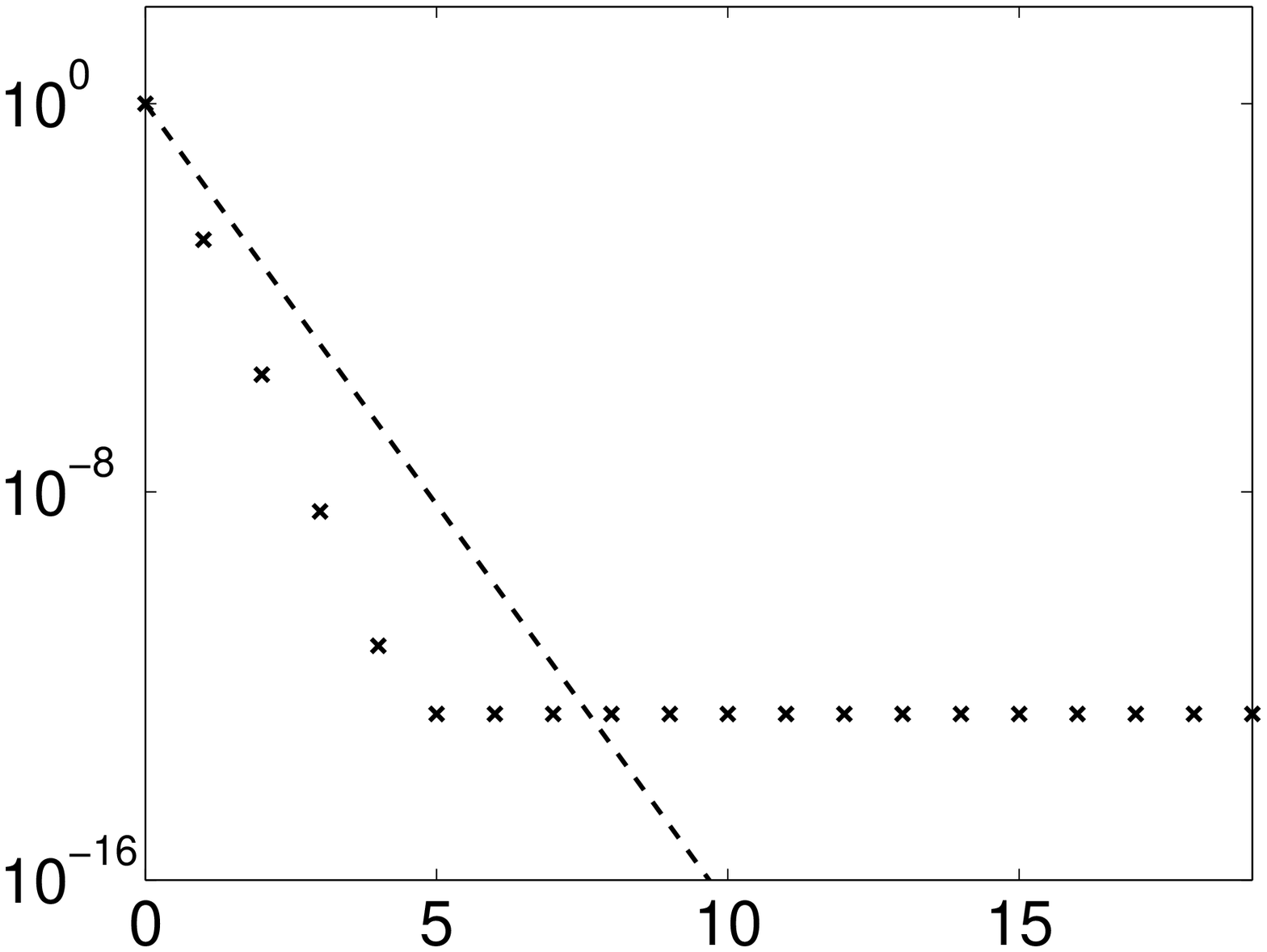}}\hfill
  \subfigure{\includegraphics[width=6cm]{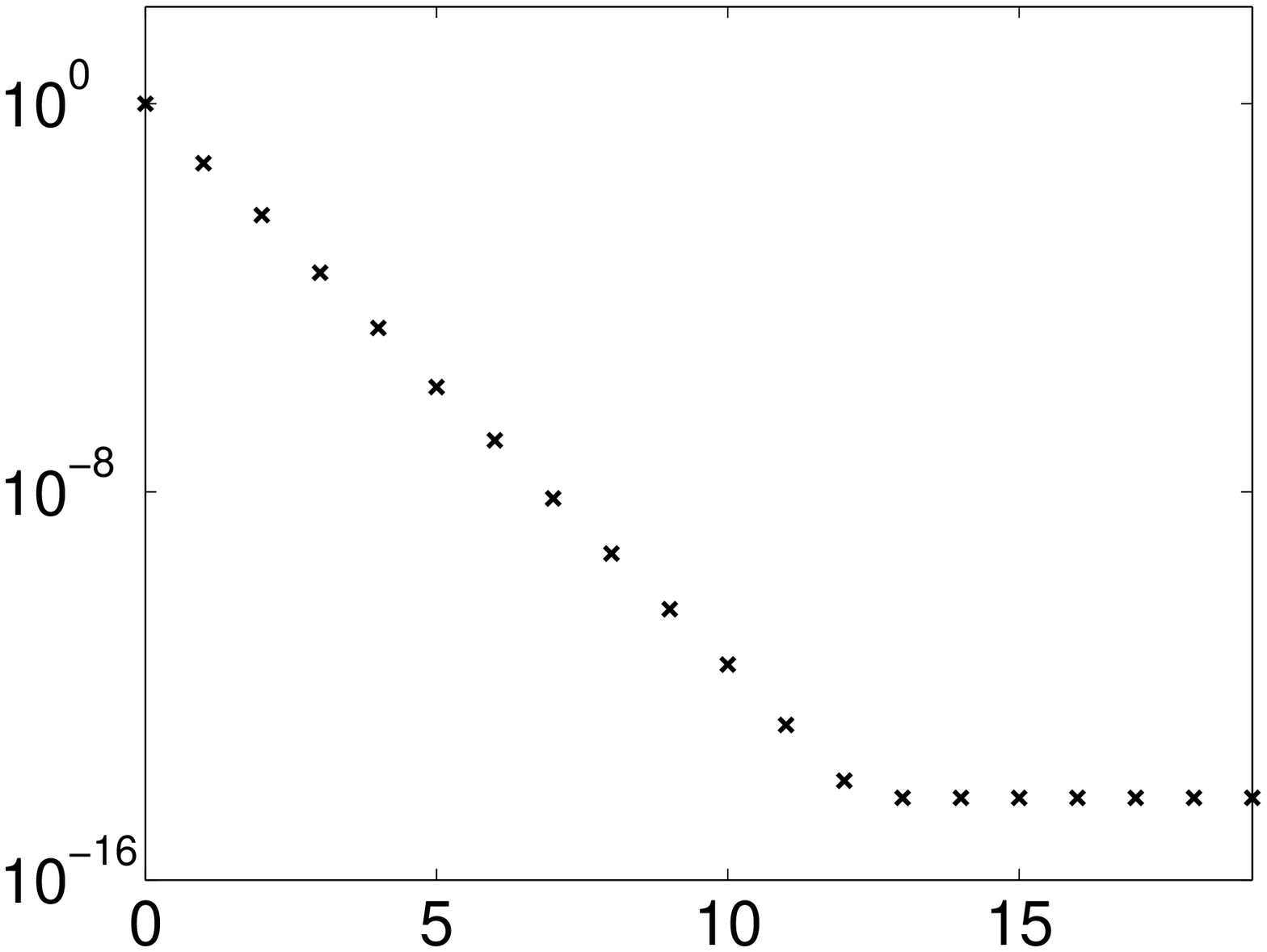}}
  \caption{Native error $\|\zb {\hat f}_l-\zb {\hat W} \zb A^{\adj} \zb
    K_N^{-1} \zb y\|_{\zb {\hat W}^{-1}}$ for the univariate interpolation
    problem with respect to the current iteration $l$.
    The number of samples is $M=100$, the number of computed Fourier
    coefficients is $N=1000$, and the separation distance of the nodes is
    $q=4\times 10^{-3}$.
    Top left: no weights, i.e., Dirichlet kernel;
    Top right: weight function $g_2$, i.e., Fej\'er kernel, predicted decay rate
    (dashed);
    Bottom left: weight function $g_4$, i.e., B-Spline kernel, predicted decay rate
    (dashed);
    Bottom right: weight function $g_{1,2,10^{-2}}$, i.e., Sobolev kernel.
    \label{fig:error_decay}}
\end{figure}

\begin{figure}[ht!]
  \centering
  \subfigure{\includegraphics[width=6cm]{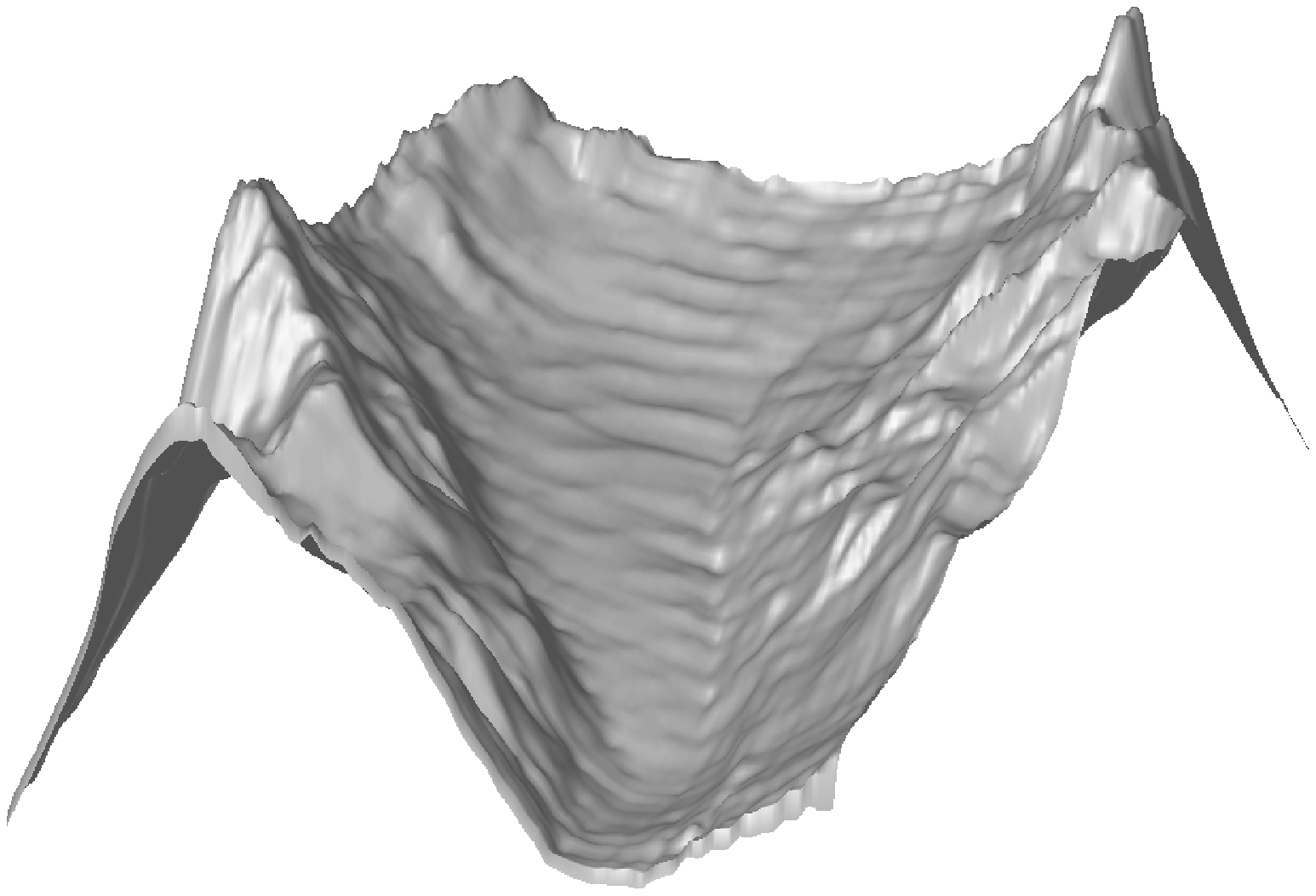}}\hfill
  \subfigure{\includegraphics[width=6cm]{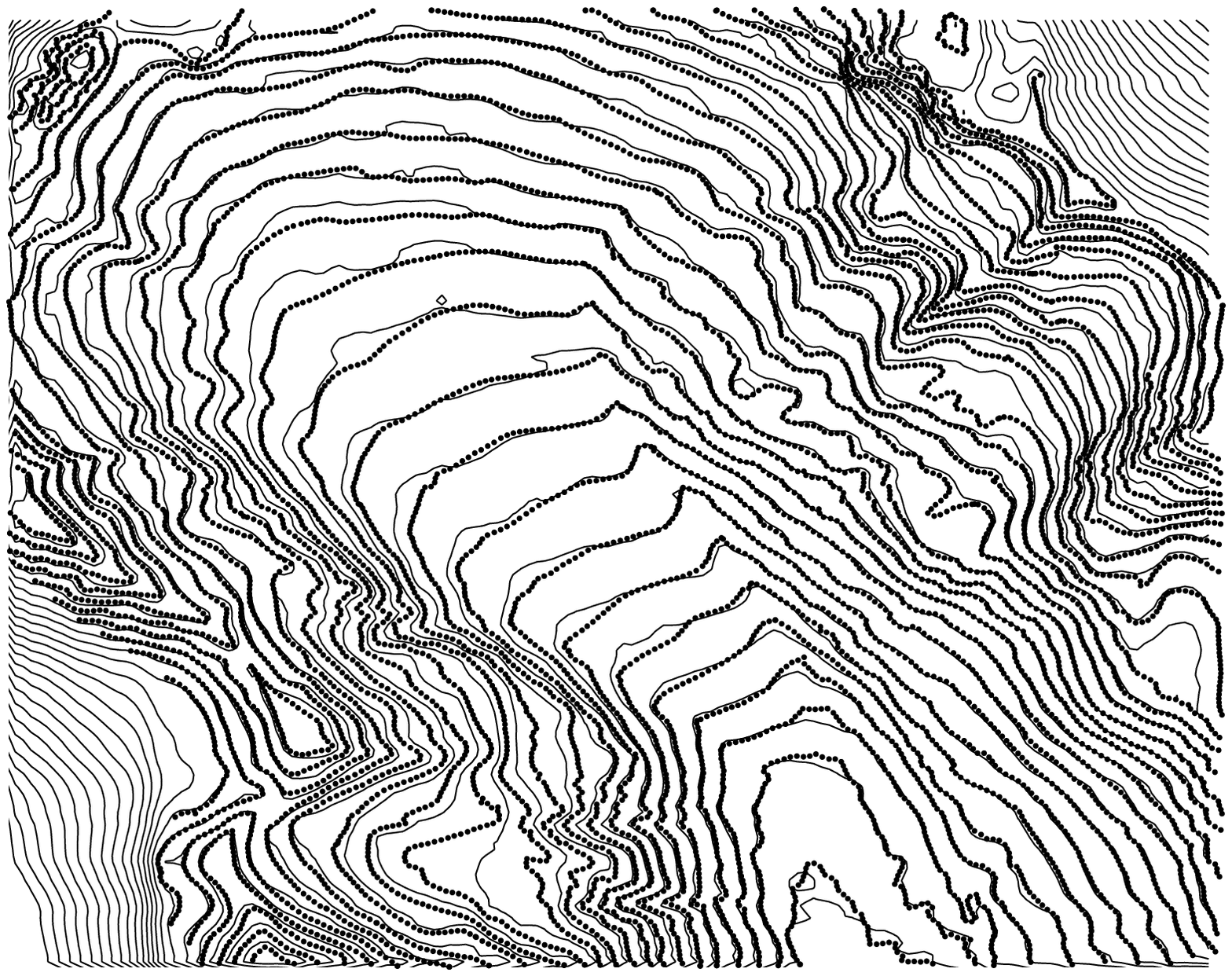}}
  \caption{Reconstruction of the glacier data set {\tt vol87.dat} from
    \cite{Franke.Daten}, $M=8345$ nodes, $N=256$, $40$ iterations, tensor
    product damping factors $\hat w_{\zb k}$ to the weight function
    $g_{\frac{1}{2},3,10^{-3}}$; see {\tt glacier\/} in \cite{kupo02C}.
    Left: surface plot, Right: contour plot and sampling set
    ($\cdot$).\label{Fig:glacier}}
\end{figure}

Furthermore, we apply Algorithm \ref{algo:CGNE} using the NFFT software
package \cite{kupo02C} to reconstruct a univariate signal from randomly
scattered data in Figure \ref{fig:error_decay} and show in Figure
\ref{Fig:glacier} the reconstruction of a bivariate signal from a
glacier data set \cite{Franke.Daten}. 
The main tool in our iterative algorithms is the NFFT, i.e., the fast matrix
times vector multiplication with $\zb A$ and $\zb A^{\adj}$, respectively.
Details concerning NFFT algorithms can be found for example in \cite{postta01}
and a software package can be found in \cite{kupo02C}.

The reconstruction of the randomly sampled univariate signal shows the decay
rates of our iterative scheme.
The sampling set consists of $M=100$ nodes separated by $q=4\times 10^{-3}$ and we
reconstruct with a polynomial degree $N=1000$ and the Dirichlet-, Fej\'er-,
B-Spline-, and Sobolev kernel.
All schemes converge within $15$ iteration where this is justified only for
the Fej\'er- and the B-Spline kernel.

The last example shows a typical test case known in  radial basis
function methods. 
We reconstruct from a data set of $M=8345$ samples on level curves of a
glacier a total number of $2^8\times 2^8 \approx 8 M$ Fourier coefficients.
Note however, that the sampling set is highly nonuniform in the sense that the
separation distance is very small compared to the mesh norm.
The assumptions of Theorem \ref{theorem:Tdallgemein} are not fulfilled.
Nevertheless, the proposed method yields a very good approximation to the
given data after $40$ iterations, which is also supported by the cross
validation test in Table \ref{tab:cv}.
Here, we exclude $\tilde M$ randomly chosen samples $\tilde {\cal
  X}\subset{\cal X}$, $\zb y_{\tilde {\cal X}}\in\C^{\tilde M}$ from the
reconstruction process and compare our approximations at these left out
nodes.
The comparison is done by means of the relative data residual and the relative
validation residual after $40$ iterations
\begin{equation*}
  r:=\frac{\|\zb y_{{\cal X}\setminus\tilde {\cal X}} -\zb A_{{\cal
  X}\setminus\tilde {\cal X}} \zb {\hat f}_{40}\|_2}{\|\zb y\|_2},\qquad
  \tilde r:=\frac{\|\zb y_{\tilde {\cal X}} -\zb A_{\tilde {\cal X}} \zb {\hat
  f}_{40}\|_2}{\|\zb y\|_2}.
\end{equation*}
As can be readily seen, the CGNE scheme achieves both a small data residual
$r$ and a small validation residual $\tilde r$.
The proposed CGNE method combines the good data fit of the CGNR scheme ($N=256$)
with the smooth approximation of the CGNR scheme ($N=64$).

\begin{table}[ht!]
  \centering
  \begin{tabular}{|c|c|c|c|c|c|c|}
    \hline
    &\multicolumn{2}{|c|}{CGNE\rule{0mm}{2.5ex}}&\multicolumn{2}{|c|}{CGNR
    ($N=256$)}&\multicolumn{2}{|c|}{CGNR ($N=64$)}\\ \hline
    $\tilde M$\rule{0mm}{3ex} & $r$ & $\tilde r$ & $r$ & $\tilde r$ & $r$ & $\tilde r$\\
    \hline
    \input{images/output_data_cv.tex}
    \hline
  \end{tabular}
  \caption{Cross validation of the reconstructions. The parameters of the CGNE
    scheme are as before. Moreover, we show the residuals of the CGNR scheme
    for the least squares problem \eqref{eq:ap} with $N=256$ (underdetermined)
    and $N=64$ (overdetermined).}\label{tab:cv}
\end{table}

\section{Conclusion}
We have shown that the optimal trigonometric interpolation problem at
$q$-separated nodes in $d$ dimensions is well conditioned for a polynomial
degree $N > 2d q^{-1}$.
However, in our further extensive numerical examples we observe that
for $N\sim M^{1/d}$ one can expect fast convergence of Algorithm
\ref{algo:CGNE}. 
If we assume furthermore a uniformity condition $q=c M^{-\frac{1}{d}}$ for the
sampling set ${\cal X}$ of cardinality $M$, then the total
arithmetical complexity for 
solving the interpolation problem \eqref{eq:optInterp} up to a prescribed
error is of order ${\cal O}\left(M \log M \right)$.

We remark that dependent on the application, one solves the weighted
approximation problem \eqref{eq:ap} or the optimal interpolation problem 
\eqref{eq:optInterp}.
Under some further mild conditions on the sampling set, both problems are
solved efficiently by means of the conjugate gradient method in conjunction
with the nonequispaced FFT.

\vspace{1cm}
{\bf Acknowledgement.}
The first author is grateful for partial support of this work by the German
Academic Exchange Service (DAAD) and the warm hospitality during his stay at
the Numerical Harmonic Analysis Group, University of Vienna.
We would also like to thank the referees for their valuable suggestions.

\bibliographystyle{abbrv}

\end{document}

%% file: images/figureTd.tex
\psset{unit=0.5cm}

\begin{pspicture}(10,10)
 \put(0.2,0.2){\framebox(9.6,9.6){}}
 \put(0.3,0.5){$R_{q,\left\lfloor q^{-1}/2 \right\rfloor}$}

 \put(1,1){\framebox(8,8){}}
 \put(1.1,1.3){$R_{q,\left\lfloor q^{-1}/2 \right\rfloor-1}$}
 \put(2,2){\framebox(6,6){}}

 \put(2.2,5.1){$\hdots$}
 \put(7.2,5.1){$\hdots$}

 \put(3,3){\framebox(4,4){}}
 \put(3.1,3.3){$R_{q,1}$}
 \put(4,4){\framebox(2,2){}}
 \put(4.1,4.3){$R_{q,0}$}

 \put(4.8,5){\line(1,0){0.4}}
 \put(5,4.8){\line(0,1){0.4}}
 \put(5.1,5.1){$\zb x_0$}

\end{pspicture}


%% file: images/figureTd_refined.tex
\psset{unit=0.5cm}

\begin{pspicture}(10,10)
 
 \put(1,1){\dashbox(8,8){}}
 \put(2,2){\framebox(6,6){}}
 \put(6,1){\line(0,1){1}}
 \put(7,1){\line(0,1){1}}
 \put(8,1){\line(0,1){1}}
 \put(8,2){\line(1,0){1}}
 \put(6,1.2){\vector(1,0){0.5}}
 \put(6.8,2){\vector(0,-1){0.5}}
 \put(7,1.2){\vector(1,0){0.5}}
 \put(7.8,2){\vector(0,-1){0.5}}
 \put(8,1.2){\vector(1,0){0.5}}
 \put(8.8,2){\vector(0,1){0.5}}
 \put(8,2.8){\vector(1,0){0.5}}

 \put(2.2,5.1){$\hdots$}
 \put(7.2,5.1){$\hdots$}


 \put(4.8,5){\line(1,0){0.4}}
 \put(5,4.8){\line(0,1){0.4}}
 \put(5.1,5.1){$\zb x_0$}
\end{pspicture}

%% file: images/output_data_cv.tex
$200$ & $6.9e-04$ & $1.7e-02$ & $5.0e-04$ & $1.4e-01$ & $8.3e-03$ & $1.7e-02$  \\
$400$ & $4.7e-04$ & $2.3e-02$ & $5.0e-04$ & $2.0e-01$ & $8.3e-03$ & $2.3e-02$  \\
$600$ & $5.7e-04$ & $2.9e-02$ & $5.1e-04$ & $2.5e-01$ & $8.1e-03$ & $2.9e-02$  \\
$800$ & $4.7e-04$ & $3.4e-02$ & $5.0e-04$ & $2.8e-01$ & $8.0e-03$ & $3.4e-02$  \\
$1000$ & $4.6e-04$ & $3.8e-02$ & $4.7e-04$ & $3.2e-01$ & $8.0e-03$ & $3.8e-02$  \\

%% file: stability.bbl
\begin{thebibliography}{10}

\bibitem{Alekseev}
V.~G. Alekseev.
\newblock {Jackson-- and Jackson--Vallee Poussin-type kernels and their
  probability applications}.
\newblock {\em Theory Probab. Appl.}, 41:137 -- 143, 1993.

\bibitem{BaGr03}
R.~F. Bass and K.~Gr{{\"o}}chenig.
\newblock Random sampling of multivariate trigonometric polynomials.
\newblock {\em {\rm SIAM} J. Math. Anal.}, 36:773 -- 795, 2004.

\bibitem{Bj96}
{\AA}.~Bj{{\"o}}rck.
\newblock {\em Numerical Methods for Least Squares Problems}.
\newblock SIAM, Philadelphia, 1996.

\bibitem{Chui88}
C.~K. Chui.
\newblock {\em Multivariate Splines}.
\newblock SIAM, Philadelphia, 1988.

\bibitem{DaDeDeM04}
I.~Daubechies, M.~Defrise, and C.~D. Mol.
\newblock An iterative thresholding algorithm for linear inverse problems with
  a sparsity constraint.
\newblock {\em Comm. Pure Appl. Math.}, 57:1413 -- 1457, 2004.

\bibitem{fas97}
H.~Fa{\ss}bender.
\newblock On numerical methods for discrete least-squares approximation by
  trigonometric polynomials.
\newblock {\em Math. Comput.}, 66:719 -- 741, 1997.

\bibitem{FeGrSt95}
H.~G. Feichtinger, K.~Gr{{\"o}}chenig, and T.~Strohmer.
\newblock Efficient numerical methods in non-uniform sampling theory.
\newblock {\em Numer. Math.}, 69:423 -- 440, 1995.

\bibitem{Franke.Daten}
R.~Franke.
\newblock http://www.math.nps.navy.mil/$\sim$rfranke/README.

\bibitem{GrSt04}
D.~Grishin and T.~Strohmer.
\newblock Fast multi-dimensional scattered data approximation with {N}eumann
  boundary conditions.
\newblock {\em Lin. Alg. Appl.}, 391:99 -- 123, 2004.

\bibitem{Groechenig92}
K.~Gr{{\"o}}chenig.
\newblock Reconstruction algorithms in irregular sampling.
\newblock {\em Math. Comput.}, 59:181 -- 194, 1992.

\bibitem{hojo}
R.~A. Horn and C.~R. Johnson.
\newblock {\em Matrix Analysis}.
\newblock Cambridge University Press, Cambridge, 1985.

\bibitem{kupo02C}
S.~Kunis and D.~Potts.
\newblock {NFFT, Softwarepackage, C subroutine library}.
\newblock http://www.tu-chemnitz.de/$\sim$potts/nfft, 2002 -- 2006.

\bibitem{KuRa06}
S.~Kunis and H.~Rauhut.
\newblock Random sampling of sparse trigonometric polynomials {II} - orthogonal
  matching pursuit versus basis pursuit.
\newblock Preprint 06-06, TU-Chemnitz, 2006.

\bibitem{MhPr00}
H.~N. Mhaskar and J.~Prestin.
\newblock On the detection of singularities of a periodic function.
\newblock {\em Adv. Comput. Math.}, 12:95 -- 131, 2000.

\bibitem{NaSiWa98}
F.~J. Narcowich, N.~Sivakumar, and J.~D. Ward.
\newblock Stability results for scattered-data interpolation on euclidean
  spheres.
\newblock {\em Adv. Comput. Math.}, 8:137 -- 163, 1998.

\bibitem{postta01}
D.~Potts, G.~Steidl, and M.~Tasche.
\newblock {Fast Fourier transforms for nonequispaced data: A tutorial}.
\newblock In J.~J. Benedetto and P.~J. S.~G. Ferreira, editors, {\em Modern
  Sampling Theory: Mathematics and Applications}, pages 247 -- 270.
  Birkh{\"a}user, Boston, 2001.

\bibitem{RaSt98}
M.~Rauth and T.~Strohmer.
\newblock Smooth approximation of potential fields from noisy scattered data.
\newblock {\em Geophysics}, 63:85 -- 94, 1998.

\bibitem{ReAmGr}
L.~Reichel, G.~S. Ammar, and W.~B. Gragg.
\newblock Discrete least squares approximation by trigonometric polynomials.
\newblock {\em Math. Comput.}, 57:273 -- 289, 1991.

\bibitem{st97}
G.~Steidl.
\newblock A note on fast {F}ourier transforms for nonequispaced grids.
\newblock {\em Adv. Comput. Math.}, 9:337 -- 353, 1998.

\bibitem{Wen05}
H.~Wendland.
\newblock {\em Scattered Data Approximation}.
\newblock Cambridge Monographs on Applied and Computational Mathematics.
  Cambridge University Press, Cambridge, 2005.

\end{thebibliography}
